\documentclass[11pt]{article}
\usepackage{mathrsfs}
\usepackage{amsfonts}
\usepackage{amsmath,amssymb}
\newcommand{\epsln}{\varepsilon}
\newcommand{\supp}{{\rm supp}}
\newcommand{\diva}{{\rm div}}

\newcommand{\rmd}{\mbox{\rm d}}
\newcommand{\bfi}{{\mbox{\boldmath $i$}}}
\newcommand{\bfu}{{\mbox{\boldmath $u$}}}
\newcommand{\bfv}{{\mbox{\boldmath $v$}}}
\newcommand{\bfw}{{\mbox{\boldmath $w$}}}

\newcommand{\rme}{{\rm e}}

\def\bgeq{\begin{equation}}
\def\edeq{\end{equation}}
\def\bgar{\begin{array}}
\def\edar{\end{array}}

\title{Sharp well-posedness and ill-posedness of the Navier-Stokes initial value problem in Besov-type spaces\thanks{This work is supported by the
  China National Natural Science Foundation under the grant number 11571381.}}
\author{{Shangbin Cui\thanks{E-mail:  cuishb@mail.sysu.edu.cn}}\\
{\small School of Mathematics, Sun Yat-Sen University, Guangzhou, Guandong 510275,}\\ [-0.2cm]
{\small People's Republic of China}\\
}
\date{}

\begin{document}
\maketitle

\begin{abstract}
  Let $B^{s,\sigma}_{pq}(\mathbb{R}^n)$ ($s\in\mathbb{R}$, $\sigma\geqslant 0$, $p,q\in [1,\infty]$) be the logrithmically refined Besov space, which
  is defined by replacing $2^{js}$ in the definition of the Besov space $B^s_{pq}(\mathbb{R}^n)$ with $2^{js}j^{\sigma}$ for all $j\in\mathbb{N}$.
  Let $B^{s,\sigma}_{\infty\,q\,0}(\mathbb{R}^n)$ ($s\in\mathbb{R}$, $\sigma\geqslant 0$, $q\in [1,\infty]$) be the closure of the Schwartz space
  $S(\mathbb{R}^n)$ in $B^{s,\sigma}_{\infty q}(\mathbb{R}^n)$. We prove that the Navier-Stokes initial value problem is locally well-posed in
  $B^{-1,\sigma}_{\infty\,q\,0}(\mathbb{R}^n)$ for $1\leqslant q\leqslant\infty$ and $\sigma\geqslant\sigma_q:=1-\min\{1-\frac{1}{q},\frac{1}{q}\}$,
  and ill-posed in $B^{-1,\sigma}_{\infty\,q}(\mathbb{R}^n)$ for $1\leqslant q\leqslant\infty$ and $0\leqslant\sigma<\sigma_q$. The well-posedness
  result is proved by using some sharp bilinear estimates obtained from some Hardy-Littlewood type inequalities. The ill-posedness assertion is
  proved by refining the arguments of Wang \cite{Wang} and Yoneda \cite{Yon10}.
\medskip

\medskip

\end{abstract}

\section{Introduction}

\hskip 2em
  This paper addresses the following question which has attracted much attention during the past two decades (cf. also Wang \cite{Wang}):
  What is the largest Besov-type space in which the initial value problem of the Navier-Stokes equations is well-posed? The purpose of this
  paper is to give an answer to this question
\footnotemark[1].
\footnotetext[1]{This is the second version of the paper with the same title publicized in ArXiv under the number 1505.00865. It remedies
  the incorrect proof of Lemma 4.2 of the previous version and some other small mistakes. The main result of this paper has been
  written in the book of Lemari\'{e}-Rieusset \cite{LEM16} as Theorem 9.6 (without proof). The author is glad to acknowledge his sincere 
  thanks to Weipeng Zhu for helping him finding the mistakes in the previous version.}

  Recall that the initial value problem of the Navier-Stokes equations reads as follows:
\begin{eqnarray}
\left\{
\begin{array}{l}
  \partial_t\bfu-\Delta\bfu+(\bfu\cdot\nabla)\bfu+\nabla\pi=0\quad \mbox{in}\;\,\mathbb{R}^n\times\mathbb{R}_+,\\
  \nabla\cdot\bfu=0\quad \mbox{in}\;\,\mathbb{R}^n\times\mathbb{R}_+,\\
  \bfu(x,0)=\bfu_0(x)\quad \mbox{for}\;\, x\in\mathbb{R}^n,
\end{array}
\right.
\end{eqnarray}
  where $n\geq 2$, $\bfu=\bfu(x,t)=(u_1(x,t),u_2(x,t),\cdots,u_n(x,t))$ is an unknown $n$-vector function in $(x,t)$ variables, $x\in\mathbb{R}^n$,
  $t\geq 0$, $\pi=\pi(x,t)$ is an unknown scalar function, $\bfu_0=\bfu_0(x)$ is a given $n$-vector function, $\Delta$ is the Laplacian in the $x$
  variables, $\nabla=(\partial_{x_1},\partial_{x_2},\cdots,\partial_{x_n})$, and $\mathbb{R}_+=(0,\infty)$.

  Let $\mathbb{P}=I+\nabla(-\Delta)^{-1}\nabla$ be the Helmholtz-Weyl projection operator, i.e., the $n\times n$ matrix pseudo-differential operator
  in $\mathbb{R}^n$ with the matrix symbol $\Big(\delta_{ij}-\frac{\xi_i\xi_j}{|\xi|^2}\Big)_{i,j=1}^n$, where $\delta_{ij}$'s are the Kronecker symbols.
  It is well-known that when only the $L^2_{uloc,x}L^2_t$-class solutions (see \cite{LEM02} for this notion) are considered, which is the case in this
  paper, the problem (1.1) is equivalent to the following formally simpler problem:
\begin{eqnarray}
\left\{
\begin{array}{l}
  \partial_t\bfu-\Delta\bfu+\mathbb{P}\nabla\cdot(\bfu\otimes\bfu)=0\quad \mbox{in}\;\,\mathbb{R}^n\times\mathbb{R}_+,\\
  \bfu(x,t)=\bfu_0(x)\quad \mbox{for}\;\, x\in\mathbb{R}^n.
\end{array}
\right.
\end{eqnarray}
  Throughout this paper, for any scaler function space $\mathscr{X}$ we shall use the same notation $\mathscr{X}$ to denote its $n$-vector counterpart
  to simplify the notation. Let $X$ be a function space continuously embedded in $S'(\mathbb{R}^n)$, the space of temperate distributions on
  $\mathbb{R}^n$ endowed with the dual topology of the Schwartz space $S(\mathbb{R}^n)$. Recall that the initial value problem $(1.1)$ is said to be
  {\em locally well-posed in $X$} if for any $\bfu_0\in X$ with $\diva\bfu_0=0$ there exists corresponding $T>0$ and a continuously embedded subspace
  $Y_T$ of $C([0,T],X)$ such that the problem $(1.2)$ has a unique solution $\bfu$ in $Y_T$, and the solution map $\bfu_0\mapsto\bfu$ is continuous
  with respect to the norm topologies of $X$ and $C([0,T],X)$. If $(1.1)$ is not locally well-posed in a function space $X$, then it is called {\em
  ill-posed in $X$}. Also recall that $(1.1)$ is said to be {\em semi-globally well-posed in $X$ for small initial data} if for any $T>0$ there exists
  corresponding constant $\varepsilon>0$ and a continuously embedded subspace $Y_T$ of $C([0,T],X)$ such that for any $\bfu_0\in X$ with $\diva\bfu_0=0$
  and $\|\bfu_0\|_X<\varepsilon$ the problem $(1.2)$ has a unique solution $\bfu$ in $Y_T$, and the solution map $\bfu_0\mapsto\bfu$ is continuous with
  respect to the norm topologies of $X$ and $C([0,T],X)$. If there exists constant $\varepsilon>0$ such that for any $\bfu_0\in X$ with $\diva\bfu_0=0$
  and $\|\bfu_0\|_X<\varepsilon$ the problem $(1.2)$ has a unique solution $\bfu$ in some subspace of $C([0,\infty),X)\cap L^{\infty}((0,\infty),X)$,
  and the solution map $\bfu_0\mapsto\bfu$ is continuous with respect to the norm topologies of $X$ and $L^{\infty}((0,\infty),X)$, then $(1.1)$ is
  said to be {\em globally well-posed in $X$ for small initial data}.

  The topic of well-posedness of the problem (1.1) in various function spaces has been deeply investigated during the past 50 years. In 1964 Fujita and
  Kato \cite{FUJK64} obtained the first result on this topic by proving that the problem $(1.1)$ is locally well-posed in $H^s(\mathbb{R}^n)$ for
  $s\geq\frac{n}{2}-1$ and globally well-posed in $H^{\frac{n}{2}-1}(\mathbb{R}^n)$ for small initial data. These results were later extended to various
  other function spaces, cf. \cite{BAR96, CAN97, FABJR72, GIG86, GM89, KAT84, KOCT01, LEM07, PLA96, TER99, WEI81} and references cited therein. Note
  that the literatures listed here are far from being complete; we refer the reader to see \cite{CAN04} and \cite{LEM02} for expositions and more
  references. Here we particularly mention that by Cannone \cite{CAN97} and Planchon \cite{PLA96}, the problem (1.1) is well-posed in the Besov spaces
  $B^{s}_{pq}(\mathbb{R}^n)$ for $s\geqslant-1+\frac{n}{p}$, $1\leqslant p<\infty$, $1\leqslant q\leqslant\infty$ (for $q=\infty$, this means that
  it is well-posed in the closure of $S(\mathbb{R}^n)$ in $B^{s}_{p\infty}(\mathbb{R}^n)$, and in what follows, similar remark should be made in
  any case where the $\infty$ index appears), and by Koch and Tataru \cite{KOCT01},
  it is well-posed in $BMO^{-1}$. Note that the inhomogeneous version $bmo^{-1}$ of $BMO^{-1}$ is the largest initial value space in which the problem
  (1.1) is known to be locally well-posed.

  On the other hand, in 2008 Bourgain and Pavlovi\'{c} \cite{BP08} proved that the problem (1.1) is ill-posed in the Besov space
  $\dot{B}^{-1}_{\infty\infty}(\mathbb{R}^n)$. Yoneda \cite{Yon10} further proved that (1.1) is ill-posed in the Besov spaces
  $\dot{B}^{-1}_{\infty q}(\mathbb{R}^n)$ and the Triebel-Lizorkin spaces $\dot{F}^{-1}_{\infty q}(\mathbb{R}^n)$ for $2<q\leqslant\infty$. Recently,
  Wang \cite{Wang} proved that the problem (1.1) is also ill-posed in the Besov spaces $\dot{B}^{-1}_{\infty q}(\mathbb{R}^n)$ for $1\leqslant
  q\leqslant 2$, which is a remarkable result because previously it had been commonly conjectured that (1.1) is well-posed in
  $\dot{B}^{-1}_{\infty q}(\mathbb{R}^n)$ for $1\leqslant q\leqslant 2$ due to the fact that they are smaller than $BMO^{-1}$. Note that all the
  above-mentioned ill-posedness results also hold for the corresponding inhomogeneous spaces, because all the arguments used in \cite{BP08},
  \cite{Yon10} and \cite{Wang} also work for the corresponding inhomogeneous spaces.

  Recalling that $BMO^{-1}=\dot{F}^{-1}_{\infty 2}(\mathbb{R}^n)$ and $bmo^{-1}=F^{-1}_{\infty 2}(\mathbb{R}^n)$, we see that  $BMO^{-1}$ and $bmo^{-1}$
  are respectively the largest homogeneous and inhomogeneous Triebel-Lizorkin spaces in which the problem (1.1) is well-posed. Naturally, we want to
  know what is the largest Besov-type space in which the problem (1.1) is well-posed. To give an answer to this question we need to refine the
  classification of the Besov space and introduce the logarithmically refined Besov space $B^{s,\sigma}_{pq}(\mathbb{R}^n)$ as
  follows (cf. \cite{Yon10}):
\medskip

  {\bf Definition 1.1} \ \ {\em $(1)$\ Let $s\in\mathbb{R}$, $\sigma\geqslant 0$ and $p,q\in[1,\infty]$. The function space
  $B^{s,\sigma}_{pq}(\mathbb{R}^n)$ consists of all temperate distributions $u$ on $\mathbb{R}^n$ such that $S_0u\in L^p(\mathbb{R}^n)$, $\Delta_j u
  \in L^p(\mathbb{R}^n)$, $j=1,2,\cdots$, and $\{2^{js}j^{-\sigma}\|\Delta_j u\|_p\}_{j=1}^\infty\in l^q$, where $S_0$ and $\Delta_j$ are the
  frequency-localizing operators appearing in the  Littlewood-Paley decomposition $u=S_0u+\sum_{j=1}^{\infty}\Delta_j u$ $($see the next section$)$.
  The norm of $u\in B^{s,\sigma}_{pq}(\mathbb{R}^n)$ is given by $\|u\|_{B^{s,\sigma}_{pq}}=\|S_0u\|_p+\|\{2^{js}j^{-\sigma}
  \|\Delta_j u\|_p\}_{j=1}^\infty\|_{l^q}$, i.e.,
$$
  \|u\|_{B^{s,\sigma}_{pq}}=\left\{
\begin{array}{ll}
  \displaystyle\|S_0u\|_p+\Big[\sum_{j=1}^\infty\Big(2^{js}j^{\sigma}\|\Delta_j u\|_p
  \Big)^q\Big]^{\frac{1}{q}} \quad &\mbox{for}\;\; 1\leqslant q<\infty,\\ [0.3cm]
  \displaystyle\|S_0u\|_p+\sup_{j\in\mathbb{N}}\Big(2^{js}j^{\sigma}\|\Delta_j u\|_p
  \Big) \quad &\mbox{for}\;\; q=\infty.
\end{array}
\right.
$$
  Here and throughout the paper $\|\cdot\|_p$ denotes the norm of $L^p(\mathbb{R}^n)$ $(1\leqslant p\leqslant\infty)$.

  $(2)$\ For $s\in\mathbb{R}$, $\sigma\geqslant 0$ and $p,q\in[1,\infty]$, we denote by $B^{s,\sigma}_{pq0}(\mathbb{R}^n)$ the closure of
  $S(\mathbb{R}^n)$ in $B^{s,\sigma}_{pq}(\mathbb{R}^n)$.}
\medskip

  It is easy to prove that $B^{s,\sigma}_{pq}(\mathbb{R}^n)$ is a Banach space, and clearly $B^{s,0}_{pq}(\mathbb{R}^n)=B^s_{pq}(\mathbb{R}^n)$, i.e.,
  when $\sigma=0$, $B^{s,\sigma}_{pq}(\mathbb{R}^n)$ coincides to the usual Besov space $B^s_{pq}(\mathbb{R}^n)$. Moreover, it is also easy to prove
  that the following embedding relations hold:
\begin{itemize}
\item For $t>s$, $\tau>\sigma>0$ and $p,q\in[1,\infty]$, we have
$$
  B^{t}_{pq}(\mathbb{R}^n)\subseteq B^{s,\tau}_{pq}(\mathbb{R}^n)\subseteq B^{s,\sigma}_{pq}(\mathbb{R}^n)\subseteq
  B^{s}_{pq}(\mathbb{R}^n)
$$
  with continuous embedding.
\item For $s\in\mathbb{R}$, $\sigma_1\geqslant\sigma_2\geqslant 0$ and $p,q_1,q_2\in[1,\infty]$ such that $\sigma_1+1/q_1>\sigma_2+1/q_2$, we have
$$
  B^{s,\sigma_1}_{pq_1}(\mathbb{R}^n)\subseteq B^{s,\sigma_2}_{pq_2}(\mathbb{R}^n)
$$
  with continuous embedding.
\end{itemize}
\medskip

  In our previous work \cite{Cui}, we proved that the problem $(1.1)$ is locally well-posed in $B^{-1,1}_{\infty\infty 0}(\mathbb{R}^n)$ and
  semi-globally well-posed in $B^{-1,1}_{\infty\infty}(\mathbb{R}^n)$ for small initial data (cf. Theorems 2.1 and 2.2 of \cite{Cui}). The arguments
  used in \cite{Cui} can be easily extended to prove that $(1.1)$ is also locally well-posed in $B^{-1,\sigma}_{\infty\infty 0}(\mathbb{R}^n)$ and
  semi-globally well-posed in $B^{-1,\sigma}_{\infty\infty}(\mathbb{R}^n)$ for small initial data for any $\sigma\geqslant 1$. Our first main result
  of this paper extends these results to $B^{-1,\sigma}_{\infty\,q\,0}(\mathbb{R}^n)$ and $B^{-1,\sigma}_{\infty\,q}(\mathbb{R}^n)$ for $1\leqslant q<
  \infty$ and $\sigma\geqslant\sigma_q$, where
$$
  \sigma_q=1-\min\Big\{1-\frac{1}{q},\frac{1}{q}\Big\} \quad \mbox{for}\;\; 1\leqslant q\leqslant\infty,
$$
  i.e., we have the following result:
\medskip

  {\bf Theorem 1.2}\ \ {\em Let $1\leqslant q<\infty$ and assume that $\sigma\geqslant\sigma_q$. Then the following assertions hold:

  $(1)$\ The problem $(1.1)$ is locally well-posed in $B^{-1,\sigma}_{\infty\,q\,0}(\mathbb{R}^n)$. More precisely, for any $\bfu_0\in
  B^{-1,\sigma}_{\infty\,q\,0}(\mathbb{R}^n)$ with $\nabla\cdot\bfu_0=0$, there exists corresponding $T>0$ such that the problem $(1.2)$ has a unique
  mild solution in the class
\begin{equation}
\left\{
\begin{array}{rcl}
  &\bfu\in C([0,T],B^{-1,\sigma}_{\infty\,q\,0}(\mathbb{R}^n))\cap L^\infty_{\rm loc}((0,T],L^\infty(\mathbb{R}^n)), \quad
  \nabla\cdot\bfu=0,& \\
  &\displaystyle\sup_{t\in (0,T)}\sqrt{t}\Big|\!\ln\Big(\frac{t}{\rme T}\Big)\Big|^{\sigma}\|\bfu(t)\|_{\infty}<\infty, \quad
  \sqrt{t}\Big|\!\ln\Big(\frac{t}{\rme T}\Big)\Big|^{\sigma}\|\bfu(t)\|_{\infty}\in L^q\Big((0,T),\frac{\rmd t}{t}\Big),&
\end{array}
\right.
\end{equation}
  and the solution map $\bfu_0\mapsto \bfu$ from a neighborhood of $\bfu_0$ in $B^{-1,\sigma}_{\infty\,q\,0}(\mathbb{R}^n)$ to the Banach space of the
  above class of functions on $\mathbb{R}^n\times (0,T)$ is Lipschitz continuous.

  $(2)$\ The problem $(1.1)$ is semi-globally well-posed in $B^{-1,\sigma}_{\infty\, q}(\mathbb{R}^n)$ for small initial data. More precisely, for any
  $T>0$ there exists corresponding constant $\varepsilon>0$ such that for any $\bfu_0\in B^{-1,\sigma}_{\infty\, q}(\mathbb{R}^n)$ with $\nabla\cdot
  \bfu_0=0$ and $\|\bfu_0\|_{B^{-1,\sigma}_{\infty\, q}}<\varepsilon$, the problem $(1.2)$ has a unique mild solution in the class
\begin{equation}
\left\{
\begin{array}{rcl}
  &\bfu\in L^\infty((0,T),B^{-1,\sigma}_{\infty\, q}(\mathbb{R}^n))\cap L^\infty_{\rm loc}((0,T],L^\infty(\mathbb{R}^n)), \quad
  \nabla\cdot\bfu=0,& \\
  &\displaystyle\sup_{t\in (0,T)}\sqrt{t}\Big|\!\ln\Big(\frac{t}{\rme T}\Big)\Big|^{\sigma}\|\bfu(t)\|_{\infty}<\infty, \quad
  \sqrt{t}\Big|\!\ln\Big(\frac{t}{\rme T}\Big)\Big|^{\sigma}\|\bfu(t)\|_{\infty}\in L^q\Big((0,T),\frac{\rmd t}{t}\Big),& \\
  &\mbox{the map $t\mapsto\bfu(t)$ is continuous with respect to $S'(\mathbb{R}^n)$-weak topology for $0<t<T$},&
\end{array}
\right.
\end{equation}
  and the solution map $\bfu_0\mapsto \bfu$ from a neighborhood of $\bfu_0$ in $B^{-1,\sigma}_{\infty\, q}(\mathbb{R}^n)$ to the Banach space of the
  above class of functions on $\mathbb{R}^n\times (0,T)$ is Lipschitz continuous.}
\medskip

  In contrast to the above result, for the case $0\leqslant\sigma<\sigma_q$ we have the following result:
\medskip

  {\bf Theorem 1.3}\ \ {\em For $1\leqslant q\leqslant\infty$ and $0\leqslant\sigma<\sigma_q$, the problem $(1.1)$ is ill-posed in
  $B^{-1,\sigma}_{\infty\, q}(\mathbb{R}^n)$. More precisely, for $0<\delta\ll1$ and $N\gg1$ there exists $\bfu_0\in S(\mathbb{R}^n)$ with
  $\|\bfu_0\|_{B^{-1,\sigma}_{\infty\,q}}\lesssim 1$ such that if we denote by $\bfu=\bfu(\delta,t)$ the solution of the problem $(1.1)$ with initial
  data $\delta\bfu_0$ $($in case such a solution exists$)$, then
\begin{equation*}
  \|\bfu(\delta,t)\|_{B^{-1,\sigma}_{\infty\,q}}\gtrsim (\ln N)^{\sigma_q-\sigma}
\end{equation*}
   for some $0<t\leqslant1/N$.}
\medskip

  From the above result and Theorems 2.1 and 2.2 of \cite{Cui} we see that the largest Besov-type spaces in which the initial value problem of the
  Navier-Stokes equations is well-posed are the spaces $B^{-1,1-\frac{1}{q}}_{\infty\,q}(\mathbb{R}^n)$, $2\leqslant q\leqslant\infty$. This answers
  the question mentioned in the beginning of this paper.

  The organization of the rest part is as follows. In the next section we make some preliminary preparations. Section 3 is devoted to giving the
  proofs of Theorem 1.2. The proof of Theorem 1.3 will be given in the last section.

\section{Preliminary preparations}
\setcounter{equation}{0}

\hskip 2em
  In this section we make some preliminary preparations.

  Choose and fix a nonnegative non-increasing function $\phi\in C^{\infty}[0,\infty)$ such that
$$
  0\leqslant\phi\leqslant1, \quad \phi(t)=1 \;\; \mbox{for}\;\; 0\leqslant t\leqslant\frac{5}{4} \quad \mbox{and} \quad
  \phi(t)=0 \;\; \mbox{for}\;\; t\geqslant\frac{3}{2},
$$
  and set
$$
  \varphi(\xi)=\phi(|\xi|), \quad \psi(\xi)=\phi(|\xi|)-\phi(2|\xi|), \quad \psi_j(\xi)=\psi(2^{-j}\xi)\;\,(j=0,1,2,\cdots)
  \quad \mbox{for}\;\; \xi\in\mathbb{R}^n.
$$
  It is easy to see that $\varphi=1$ on $\bar{B}(0,5/4)$ and $\supp\varphi\subseteq\bar{B}(0,3/2)$, $\psi=1$ on $\bar{B}(0,5/4)\backslash B(0,3/4)$
  and $\supp\psi\subseteq\bar{B}(0,3/2)\backslash B(0,5/8)$. Here $B(a,r)$ and $\bar{B}(a,r)$ ($a\in\mathbb{R}^n$, $r>0$) respectively represent the
  open and closed balls in $\mathbb{R}^n$ with center $a$ and radius $r$. We also note that
$$
  \psi_j=1 \quad \mbox{on}\;\; \bar{B}(0,5\cdot 2^{j-2})\backslash B(0,3\cdot 2^{j-2}) \quad \mbox{and} \quad
  \supp\psi_j\subseteq B(0, 2^{j+1})\backslash\bar{B}(0, 2^{j-1})
$$
  $(j=0,1,2,\cdots)$. Moreover we have
$$
  \varphi(\xi)+\sum_{j=1}^{\infty}\psi_j(\xi)=1 \quad \mbox{for}\;\; \xi\in\mathbb{R}^n.
$$
  We denote by $\;\hat{}\;$ and $\mathscr{F}$ the Fourier transform, and by $\;\check{}\;$ and $\mathscr{F}^{-1}$ the inverse Fourier transform. The
  notation $O_M(\mathbb{R}^n)$ denotes the topological vector space of temperate smooth functions on $\mathbb{R}^n$, i.e. $u\in O_M(\mathbb{R}^n)$
  if and only if $u\in C^{\infty}(\mathbb{R}^n)$ and for any $\alpha\in\mathbb{Z}_+^n$, there exists corresponding $r\in\mathbb{R}$ such that
  $|\partial^{\,\alpha}u(x)|\lesssim (1+|x|)^r$ for $x\in\mathbb{R}^n$. Then we define $S_0:S'(\mathbb{R}^n)\to O_M(\mathbb{R}^n)$ and $\Delta_j:
  S'(\mathbb{R}^n)\to O_M(\mathbb{R}^n)$ $(j=0,1,2,\cdots)$ to be the following operators:
$$
  S_0(u)=\mathscr{F}^{-1}(\varphi\hat{u}), \quad \Delta_j(u)=\mathscr{F}^{-1}(\psi_j\hat{u}) \quad \mbox{for}\;\; u\in S'(\mathbb{R}^n)
$$
  $(j=0,1,2,\cdots)$. It is well-known that for any $u\in S'(\mathbb{R}^n)$ there holds the relation
$$
  S_0(u)+\sum_{j=1}^{\infty}\Delta_j(u)=u
$$
  in $S'(\mathbb{R}^n)$-weak topology.

  As usual for $t\geqslant 0$ we denote by $\rme^{t\Delta}$ the pseudo-differential operator on $\mathbb{R}^n$ with symbol $\rme^{-t|\xi|^2}$, i.e.,
  $\rme^{t\Delta}$ is the continuous linear operator in $S'(\mathbb{R}^n)$ defined by
$$
  \rme^{t\Delta}u=\mathscr{F}^{-1}(\rme^{-t|\xi|^2}\hat{u}(\xi)) \quad \mbox{for}\;\; u\in S'(\mathbb{R}^n).
$$
  It is well-known that when restricted on shift-invariant Banach space of test functions (see \cite{LEM02} for this concept), the family of operators
  $\{\rme^{t\Delta}\}_{t\geqslant 0}$ forms a $C_0$-semigroup of contractions (i.e. $\|\rme^{t\Delta}u\|\leqslant\|u\|$ for all $t\geqslant 0$), and
  when restricted on shift-invariant Banach space of distributions (also see \cite{LEM02} for this concept), $\{\rme^{t\Delta}\}_{t\geqslant 0}$ is a
  semigroup of contractions, but it is not necessarily strongly continuous at $t=0$ (it is strongly continuous for $t>0$).

  In the proofs of Theorems 1.3 and 1.6 we shall use the following characterization of the space $B^{s,\sigma}_{pq}(\mathbb{R}^n)$:
\medskip

  {\bf Lemma 2.1}\ \ {\em Let $s\in\mathbb{R}$, $\sigma\geqslant 0$ and $p,q\in[1,\infty]$. Let $u\in S'(\mathbb{R}^n)$. Let $t_0>0$ be
  given. Let $\gamma\geqslant 0$ and $\gamma>s$. Then $u\in B^{s,\sigma}_{pq}(\mathbb{R}^n)$ if and only if for any $t>0$ we have $e^{t\Delta}u\in
  L^p(\mathbb{R}^n)$ and $t^{-\frac{s}{2}}|\ln(\frac{t}{\rme t_0})|^{\sigma}\|(\sqrt{-t\Delta})^{\gamma}e^{t\Delta}u\|_q\in L^q((0,t_0),\frac{dt}{t})$.
  Moreover, the norms $\|u\|_{B^{s,\sigma}_{pq}}$ and
$$
  \|u\|_{B^{s,\sigma}_{pq,t_0}}=\left\{
\begin{array}{ll}
  \displaystyle\|e^{t_0\Delta}u\|_p+\Big[\int_0^{t_0}\!\Big(t^{-\frac{s}{2}}\Big|\ln\Big(\frac{t}{\rme t_0}\Big)\Big|^{\sigma}
  \|(\sqrt{-t\Delta})^{\gamma}e^{t\Delta}u\|_p\Big)^q\frac{dt}{t}\Big]^{\frac{1}{q}} \quad &\mbox{for}\;\; 1\leqslant q<\infty\\ [0.3cm]
  \displaystyle\|e^{t_0\Delta}u\|_p+
  \sup_{0<t<t_0}t^{-\frac{s}{2}}\Big|\ln\Big(\frac{t}{\rme t_0}\Big)\Big|^{\sigma}\|(\sqrt{-t\Delta})^{\gamma}e^{t\Delta}u\|_p \quad
  &\mbox{for}\;\; q=\infty
\end{array}
\right.
$$
  are equivalent.}
\medskip

  The proof is not hard; one needs only to slightly modify the arguments used in the proof of Theorem 5.3 in \cite{LEM02} to fit the present
  situation. We omit it here. $\quad\Box$

  Note that if $s<0$ and $\gamma=0$ then for $0<t_1<t_2$ we have
$$
  \|u\|_{B^{s,\sigma}_{pq,t_1}}\leqslant\|u\|_{B^{s,\sigma}_{pq,t_2}}\leqslant
  \Big(\frac{t_2}{t_1}\Big)^{\frac{|s|}{2}}\|u\|_{B^{s,\sigma}_{pq,t_1}} \quad \mbox{for}\;\; u\in B^{s,\sigma}_{pq}(\mathbb{R}^n).
$$

  It is well-known that the problem (1.2) is equivalent to the following integral equation:
$$
  \bfu(t)=\rme^{t\Delta}\bfu_0+\int_0^t\rme^{(t-\tau)\Delta}\mathbb{P}\nabla\cdot[\bfu(\tau)\otimes\bfu(\tau)]\rmd\tau.
$$
  Given $T>0$, let $B$ be the following bilinear form:
$$
  B(\bfu,\bfv)(t)=\int_0^t\rme^{(t-\tau)\Delta}\mathbb{P}\nabla\cdot[\bfu(\tau)\otimes\bfv(\tau)]\rmd\tau.
$$
  The following very useful preliminary result is well-known (cf. Chapter 15 of \cite{LEM02}):
\medskip

  {\bf Lemma 2.2}\ \ {\em For any $\bfu,\bfv\in L^1_{loc}((0,T),[L^{\infty}(\mathbb{R}^n)]^n)$ such that the right-hand side makes sense for
  a. e. $t\in (0,T)$, the following estimate holds:}
$$
  \;\;\qquad\qquad\|B(\bfu,\bfv)(t)\|_{\infty}\lesssim\int_0^t(t-\tau)^{-\frac{1}{2}}\|\bfu(\tau)\|_{\infty}\|\bfv(\tau)\|_{\infty}\rmd\tau \quad
  \mbox{for a. e. }\; t\in (0,T). \qquad\qquad\;\; \Box
$$

  To make estimate of the right-hand of the above inequality, we need some Hardy-Littlewood type inequalities, which are given in the following two
  lemmas:
\medskip

  {\bf Lemma 2.3}\ \ {\em Let $T>0$ and $F(t)=\displaystyle\int_0^t\!f(\tau)\frac{\rmd\tau}{\tau}$, $0<t<T$, where $f$ is a measurable function
  defined in $(0,T)$. Then for any $2\leqslant q\leqslant\infty$ the following estimate holds:
\begin{equation}
  \Big[\int_0^T\!\!\Big|\ln\Big(\frac{t}{\rme T}\Big)\Big|^{q-1}|F(t)|^q\frac{\rmd t}{t}\Big]^{\frac{1}{q}}\lesssim_q
  \Big[\int_0^T\!\!\Big|\ln\Big(\frac{t}{\rme T}\Big)\Big|^{q-1}|f(t)|^{\frac{q}{2}}\frac{\rmd t}{t}\Big]^{\frac{2}{q}},
\end{equation}
  where for $q=\infty$ the integration is understood in the conventional way.}
\medskip

  {\em Proof}:\ \ First, by using the Minkowsky inequality we have
\begin{eqnarray*}
  \Big[\int_0^T\!\!\Big|\ln\Big(\frac{t}{\rme T}\Big)\Big||F(t)|^2\frac{\rmd t}{t}\Big]^{\frac{1}{2}}
  &\;\leqslant\;&\Big[\int_0^T\!\!\Big(\int_0^t\Big|\ln\Big(\frac{t}{\rme T}\Big)\Big|^{\frac{1}{2}}|f(\tau)|\frac{\rmd\tau}{\tau}\Big)^2
  \frac{\rmd t}{t}\Big]^{\frac{1}{2}}
\\
  &\leqslant&\int_0^T\!\!\Big(\int_{\tau}^T\Big|\ln\Big(\frac{t}{\rme T}\Big)\Big|\frac{\rmd t}{t}\Big)^{\frac{1}{2}}|f(\tau)|\frac{\rmd\tau}{\tau}
\\
  &\lesssim&\int_0^T\!\!\Big|\ln\Big(\frac{\tau}{\rme T}\Big)\Big||f(\tau)|\frac{\rmd\tau}{\tau}.
\end{eqnarray*}
  This proves that (2.1) holds for $q=2$. Next we have
\begin{eqnarray*}
  |F(t)|&\;\leqslant\;&\int_0^t\!|f(\tau)|\frac{\rmd\tau}{\tau}=\int_0^t\!\Big|\ln\Big(\frac{\tau}{\rme T}\Big)\Big|^2|f(\tau)|\cdot
  \Big|\ln\Big(\frac{\tau}{\rme T}\Big)\Big|^{-2}\frac{\rmd\tau}{\tau}
\\
  &\leqslant&\sup_{0<\tau<T}\Big|\ln\Big(\frac{\tau}{\rme T}\Big)\Big|^2|f(\tau)|\cdot
  \int_0^t\Big|\ln\Big(\frac{\tau}{\rme T}\Big)\Big|^{-2}\,\frac{\rmd\tau}{\tau}
\\
  &\leqslant&\sup_{0<\tau<T}\Big|\ln\Big(\frac{\tau}{\rme T}\Big)\Big|^2|f(\tau)|\cdot\Big|\ln\Big(\frac{t}{\rme T}\Big)\Big|^{-1}, \quad
  \forall t\in (0,T),
\end{eqnarray*}
  so that
$$
  \sup_{0<t<T}\Big|\ln\Big(\frac{t}{\rme T}\Big)\Big||F(t)|\lesssim
  \sup_{0<t<T}\Big|\ln\Big(\frac{t}{\rme T}\Big)\Big|^2|f(t)|,
$$
  showing that (2.1) also holds for $q=\infty$. Hence, by interpolation we see that (2.1) holds for all $2\leqslant q\leqslant\infty$. $\quad\Box$
\medskip

  {\bf Lemma 2.4}\ \ {\em Let $T>0$ and $F(t)=\displaystyle\int_0^t\Big|\ln\Big(\frac{\tau}{\rme T}\Big)\Big|^{-1}\!f(\tau)\frac{\rmd\tau}{\tau}$,
  $0<t<T$, where $f$ is a measurable function defined in $(0,T)$. Then for any $1\leqslant q\leqslant\infty$ the following estimate holds:
\begin{equation}
  \Big\{\int_0^T\!\!\Big[\Big|\ln\Big(\frac{t}{\rme T}\Big)\Big||F(t)|\Big]^q\frac{\rmd t}{t}\Big\}^{\frac{1}{q}}\lesssim
  \Big\{\int_0^T\!\!\Big[\Big|\ln\Big(\frac{t}{\rme T}\Big)\Big||f(t)|\Big]^q\frac{\rmd t}{t}\Big\}^{\frac{1}{q}},
\end{equation}
  where for $q=\infty$ the integration is understood in the conventional way.}
\medskip

  {\em Proof}:\ \ First we have
\begin{eqnarray*}
  \int_0^T\!\Big|\ln\Big(\frac{t}{\rme T}\Big)\Big||F(t)|\frac{\rmd t}{t}&\;\leqslant\;&
  \int_0^T\!\!\Big[\int_0^t\Big|\ln\Big(\frac{t}{\rme T}\Big)\Big|\Big|\ln\Big(\frac{\tau}{\rme T}\Big)\Big|^{-1}\!|f(\tau)|\frac{\rmd\tau}{\tau}
  \Big]\frac{\rmd t}{t}
\\
  &=&\int_0^T\!\!\Big[\int_{\tau}^T\Big|\ln\Big(\frac{t}{\rme T}\Big)\Big|\frac{\rmd t}{t}\Big]
  \Big|\ln\Big(\frac{\tau}{\rme T}\Big)\Big|^{-1}\!|f(\tau)|\frac{\rmd\tau}{\tau}
\\
  &\leqslant&\frac{1}{2}\int_0^T\!\Big|\ln\Big(\frac{\tau}{\rme T}\Big)\Big|^2\cdot\Big|\ln\Big(\frac{\tau}{\rme T}\Big)\Big|^{-1}\!
  |f(\tau)|\frac{\rmd\tau}{\tau}
\\
  &=&\frac{1}{2}\int_0^T\!\Big|\ln\Big(\frac{\tau}{\rme T}\Big)\Big||f(\tau)|\frac{\rmd\tau}{\tau},
\end{eqnarray*}
  showing that (2.2) holds for $q=1$. Next we have
\begin{eqnarray*}
  \sup_{0<t<T}\Big|\ln\Big(\frac{t}{\rme T}\Big)\Big||F(t)|&\;\leqslant\;&
  \sup_{0<t<T}\Big|\ln\Big(\frac{t}{\rme T}\Big)\Big|\int_0^t\Big|\ln\Big(\frac{\tau}{\rme T}\Big)\Big|^{-1}\!|f(\tau)|\frac{\rmd\tau}{\tau}
\\
  &\leqslant&\sup_{0<t<T}\Big|\ln\Big(\frac{t}{\rme T}\Big)\Big|\int_0^t\Big|\ln\Big(\frac{\tau}{\rme T}\Big)\Big|^{-2}\,\frac{\rmd\tau}{\tau}
  \cdot\sup_{0<\tau<T}\Big|\ln\Big(\frac{\tau}{\rme T}\Big)\Big||f(\tau)|
\\
  &\leqslant&\sup_{0<t<T}\Big|\ln\Big(\frac{t}{\rme T}\Big)\Big||f(t)|,
\end{eqnarray*}
  showing that (2.2) also holds for $q=\infty$. Hence, by interpolation we see that (2.2) holds for all $1\leqslant q\leqslant\infty$. $\quad\Box$
\medskip

  Now let $\mathcal{B}$ be the following bilinear operator:
$$
  \mathcal{B}(f,g)(t)=\int_0^t(t-\tau)^{-\frac{1}{2}}f(\tau)g(\tau)\rmd\tau, \quad t>0.
$$
  For $1\leqslant q\leqslant\infty$, $\sigma\geqslant 0$ and $T>0$, we denote by $\dot{\mathscr{K}}^{\sigma}_q(T)$ the following function space on
  $(0,T)$:
$$
  \dot{\mathscr{K}}^{\sigma}_q(T)=\Big\{f\in L^1_{\rm loc}(0,T]:\sqrt{t}\Big|\ln\Big(\frac{t}{\rme T}\Big)\Big|^{\sigma}f(t)\in
  L^q\Big((0,T);\frac{\rmd t}{t}\Big)\Big\},
$$
  with norm
$$
  \|f\|_{\dot{\mathscr{K}}^{\sigma}_q}=\Big\|\sqrt{t}\Big|\ln\Big(\frac{t}{\rme T}\Big)\Big|^{\sigma}f(t)\Big\|_{L^q((0,T);\frac{dt}{t})},
$$
  and by $\mathscr{K}^{\sigma}_q(T)$ the following function space on $(0,T)$: $\mathscr{K}^{\sigma}_{\infty}(T)=\dot{\mathscr{K}}^{\sigma}_{\infty}(T)$
  and for $1\leqslant q<\infty$,
$$
  \mathscr{K}^{\sigma}_q(T)=\dot{\mathscr{K}}^{\sigma}_q(T)\cap\dot{\mathscr{K}}^{\sigma}_{\infty}(T),
$$
  with norm
$$
  \|f\|_{\mathscr{K}^{\sigma}_q}=\|f\|_{\dot{\mathscr{K}}^{\sigma}_q}+\|f\|_{\dot{\mathscr{K}}^{\sigma}_{\infty}}.
$$
  It is easy to prove that both $\dot{\mathscr{K}}^{\sigma}_q(T)$ and $\mathscr{K}^{\sigma}_q(T)$ ($1\leqslant q\leqslant\infty$) are Banach spaces.
  We use the notation $\mathscr{K}^{\sigma}_{\infty\,0}(T)$ to denote the closure of the subspace of simple functions in $[0,T]$ in
  $\mathscr{K}^{\sigma}_{\infty}(T)$. From the interpolation theory for $L^p$-spaces with change of measures (cf. Sections 5.4 and 5.5 of
  \cite{BerLof}) we see that the following interpolation relations hold: For $1\leqslant q_0,q_1<\infty$,
$$
  [\dot{\mathscr{K}}^{\sigma_0}_{q_0}(T),\dot{\mathscr{K}}^{\sigma_1}_{q_1}(T)]_{[\theta]}=\dot{\mathscr{K}}^{\sigma}_q(T), \qquad
  [\mathscr{K}^{\sigma_0}_{q_0}(T),\mathscr{K}^{\sigma_1}_{q_1}(T)]_{[\theta]}=\mathscr{K}^{\sigma}_q(T),
$$
  and for $1\leqslant q_0<\infty$, $q_1=\infty$,
$$
  [\dot{\mathscr{K}}^{\sigma_0}_{q_0}(T),\dot{\mathscr{K}}^{\sigma_1}_{\infty\,0}(T)]_{[\theta]}=\dot{\mathscr{K}}^{\sigma}_q(T), \qquad
  [\mathscr{K}^{\sigma_0}_{q_0}(T),\mathscr{K}^{\sigma_1}_{\infty\,0}(T)]_{[\theta]}=\mathscr{K}^{\sigma}_q(T),
$$
  where $0<\theta<1$, $\sigma=(1\!-\!\theta)\sigma_0+\theta\sigma_1$, $1/q=(1\!-\!\theta)/q_0+\theta/q_1$. Besides, it is clear that
$$
  \mathscr{K}^{\sigma}_{1}(T)\subseteq\mathscr{K}^{\sigma}_{q_1}(T)\subseteq\mathscr{K}^{\sigma}_{q_2}(T)\subseteq\mathscr{K}^{\sigma}_{\infty\,0}(T)
  \quad \mbox{for}\;\; 1\leqslant q_1\leqslant q_2<\infty, \quad \sigma\geqslant 0,
$$
  and
$$
  \dot{\mathscr{K}}^{\sigma_1}_q(T)\subseteq\dot{\mathscr{K}}^{\sigma_2}_q(T), \qquad
  \mathscr{K}^{\sigma_1}_q(T)\subseteq\mathscr{K}^{\sigma_2}_q(T)
  \quad \mbox{for}\;\; \sigma_1\geqslant\sigma_2\geqslant 0, \quad 1\leqslant q\leqslant\infty.
$$
  The following bilinear estimate will play a fundamental role in the proof of Theorem 1.2:
\medskip

  {\bf Lemma 2.5}\ \ {\em Let $T>0$ be given and assume that $(q,\sigma)$ satisfies one of the following two conditions: $(a)$ $\sigma\geqslant 1$ and
  $1\leqslant q\leqslant\infty$;  $(b)$ $1/2\leqslant\sigma<1$ and $1/\sigma\leqslant q\leqslant 1/(1\!-\!\sigma)$. Then the following estimate
  holds:}
\begin{equation}
  \|\mathcal{B}(f,g)\|_{\mathscr{K}^{\sigma}_q}\lesssim_{q,\sigma}\|f\|_{\mathscr{K}^{\sigma}_q}\|g\|_{\mathscr{K}^{\sigma}_q}, \quad
  \forall f,g\in\mathscr{K}^{\sigma}_q(T).
\end{equation}

  {\em Proof}:\ \ First we have
\begin{eqnarray*}
  |\mathcal{B}(f,g)(t)|&\;\leqslant\;&\int_0^t(t-\tau)^{-\frac{1}{2}}|f(\tau)||g(\tau)|\rmd\tau
\\
  &\leqslant&\sqrt{\frac{2}{t}}\int_0^t|f(\tau)||g(\tau)|\rmd\tau+\int_{\frac{t}{2}}^t(t-\tau)^{-\frac{1}{2}}|f(\tau)||g(\tau)|\rmd\tau
  :=J_1(t)+J_2(t).
\end{eqnarray*}
  For $J_1(t)$ we have
$$
  \sqrt{t}J_1(t)=\sqrt{2}\int_0^t\sqrt{\tau}|f(\tau)|\cdot\sqrt{\tau}|g(\tau)|\frac{\rmd\tau}{\tau}, \quad 0<t<T.
$$
  Applying (2.1) to the cases $q=\infty$ and $q=2$ we respectively get
\begin{eqnarray}
  &\|J_1\|_{\dot{\mathscr{K}}^1_{\infty}}\lesssim\|f\|_{\dot{\mathscr{K}}^1_{\infty}}\|g\|_{\dot{\mathscr{K}}^1_{\infty}}, \quad
  \forall f,g\in\dot{\mathscr{K}}^1_{\infty}(T),&
\\
  &\|J_1\|_{\dot{\mathscr{K}}^{\frac{1}{2}}_2}\lesssim\|f\|_{\dot{\mathscr{K}}^{\frac{1}{2}}_2}\|g\|_{\dot{\mathscr{K}}^{\frac{1}{2}}_2}, \quad
  \forall f,g\in\dot{\mathscr{K}}^{\frac{1}{2}}_2(T).&
\end{eqnarray}
  Next we note that
$$
  \sqrt{t}J_1(t)\leqslant\sqrt{2}\sup_{0<\tau<T}\sqrt{\tau}\Big|\ln\Big(\frac{\tau}{\rme T}\Big)\Big||f(\tau)|\cdot
  \int_0^t\Big|\ln\Big(\frac{\tau}{\rme T}\Big)\Big|^{-1}\!\sqrt{\tau}|g(\tau)|\frac{\rmd\tau}{\tau}, \quad 0<t<T,
$$
  so that by applying (2.2) to the case $q=1$ we obtain
\begin{equation}
  \|J_1\|_{\dot{\mathscr{K}}^1_1}\lesssim\|f\|_{\dot{\mathscr{K}}^1_{\infty}}\|g\|_{\dot{\mathscr{K}}^1_1}, \quad
  \forall f\in\dot{\mathscr{K}}^1_{\infty}(T),\;\; \forall g\in\dot{\mathscr{K}}^1_1(T).
\end{equation}
  By using bilinear interpolation (cf. Theorem 4.4.1 of \cite{BerLof}), from (2.4) $\sim$ (2.6) we easily get the following estimate provided that
  $1/2\leqslant\sigma\leqslant 1$ and $1/\sigma\leqslant q\leqslant 1/(1\!-\!\sigma)$:
$$
  \|J_1\|_{\dot{\mathscr{K}}^{\sigma}_q}\lesssim\|f\|_{\dot{\mathscr{K}}^{\sigma}_{r_{\sigma}}}\|g\|_{\dot{\mathscr{K}}^{\sigma}_q}, \quad
  \forall f\in\dot{\mathscr{K}}^{\sigma}_{r_{\sigma}}(T),\;\; \forall g\in\dot{\mathscr{K}}^{\sigma}_q(T),
$$
  where $r_{\sigma}=1/(1\!-\!\sigma)$, which immediately implies that
\begin{equation}
  \|J_1\|_{\mathscr{K}^{\sigma}_q}\lesssim\|f\|_{\mathscr{K}^{\sigma}_q}\|g\|_{\mathscr{K}^{\sigma}_q}, \quad
  \forall f,g\in\mathscr{K}^{\sigma}_q(T),
\end{equation}
  provided that $1/2\leqslant\sigma\leqslant 1$ and $1/\sigma\leqslant q\leqslant 1/(1\!-\!\sigma)$.  For $J_2(t)$, since $\tau\sim t$, we have
\begin{eqnarray*}
  \sup_{0<t<T}\sqrt{t}\Big|\ln\Big(\frac{t}{\rme T}\Big)\Big|^{\sigma}J_2(t)
  &\;\leqslant\;&\sup_{0<\tau<T}\sqrt{\tau}|f(\tau)|\cdot\sup_{0<\tau<T}\sqrt{\tau}\Big|\ln\Big(\frac{\tau}{\rme T}\Big)\Big|^{\sigma}|g(\tau)|\cdot
  \int_{\frac{t}{2}}^t(t-\tau)^{-\frac{1}{2}}\tau^{-\frac{1}{2}}\rmd\tau
\\
  &\;\lesssim\;&\sup_{0<\tau<T}\sqrt{\tau}|f(\tau)|\cdot\sup_{0<\tau<T}\sqrt{\tau}\Big|\ln\Big(\frac{\tau}{\rme T}\Big)\Big|^{\sigma}|g(\tau)|,
  \quad \forall\sigma\geqslant 0.
\end{eqnarray*}
  Hence
\begin{equation}
  \|J_2\|_{\dot{\mathscr{K}}^{\sigma}_{\infty}}\lesssim\|f\|_{\dot{\mathscr{K}}^{0}_{\infty}}\|g\|_{\dot{\mathscr{K}}^{\sigma}_{\infty}}
  \lesssim\|f\|_{\dot{\mathscr{K}}^{\sigma}_{\infty}}\|g\|_{\dot{\mathscr{K}}^{\sigma}_{\infty}},
  \quad \forall f,g\in\dot{\mathscr{K}}^{\sigma}_{\infty}(T), \quad \forall\sigma\geqslant 0.
\end{equation}
  Besides, since $\tau\sim t$ also implies that
$$
  \sqrt{t}\Big|\ln\Big(\frac{t}{\rme T}\Big)\Big|^{\sigma}J_2(t)\leqslant\sup_{0<\tau<T}\sqrt{\tau}|f(\tau)|\cdot
  \int_{\frac{t}{2}}^t(t-\tau)^{-\frac{1}{2}}\tau^{-\frac{1}{2}}\cdot\sqrt{\tau}\Big|\ln\Big(\frac{\tau}{\rme T}\Big)\Big|^{\sigma}|g(\tau)|\rmd\tau,
$$
  we further have
\begin{eqnarray*}
  \int_0^T\!\!\sqrt{t}\Big|\ln\Big(\frac{t}{\rme T}\Big)\Big|^{\sigma}J_2(t)\frac{\rmd t}{t}
  &\;\leqslant\;&\sup_{0<\tau<T}\sqrt{\tau}|f(\tau)|\cdot\int_0^T\!\Big(\int_{\tau}^{2\tau}\!(t-\tau)^{-\frac{1}{2}}t^{-1}\rmd t\Big)
  \tau^{-\frac{1}{2}}\cdot\sqrt{\tau}\Big|\ln\Big(\frac{\tau}{\rme T}\Big)\Big|^{\sigma}|g(\tau)|\rmd\tau
\\
  &\;\lesssim\;&\sup_{0<\tau<T}\sqrt{\tau}|f(\tau)|\cdot\int_0^T\!\!\sqrt{\tau}\Big|\ln\Big(\frac{\tau}{\rme T}\Big)\Big|^{\sigma}|g(\tau)|
  \frac{\rmd\tau}{\tau}, \quad \forall\sigma\geqslant 0,
\end{eqnarray*}
  i.e.,
\begin{equation}
  \|J_2\|_{\dot{\mathscr{K}}^{\sigma}_1}\lesssim\|f\|_{\dot{\mathscr{K}}^{0}_{\infty}}\|g\|_{\dot{\mathscr{K}}^{\sigma}_1}
  \lesssim\|f\|_{\dot{\mathscr{K}}^{\sigma}_{\infty}}\|g\|_{\dot{\mathscr{K}}^{\sigma}_1},
  \quad \forall f\in\dot{\mathscr{K}}^{\sigma}_{\infty}(T),\;\; \forall g\in\dot{\mathscr{K}}^{\sigma}_1(T),
   \quad \forall\sigma\geqslant 0.
\end{equation}
  By interpolation, from (2.8) and (2.9) we get
$$
  \|J_2\|_{\dot{\mathscr{K}}^{\sigma}_q}\lesssim\|f\|_{\dot{\mathscr{K}}^{\sigma}_{\infty}}\|g\|_{\dot{\mathscr{K}}^{\sigma}_q},
  \quad \forall f\in\dot{\mathscr{K}}^{\sigma}_{\infty}(T),\;\; \forall g\in\dot{\mathscr{K}}^{\sigma}_q(T),
   \quad \forall q\in [1,\infty],\;\; \forall \sigma\geqslant 0,
$$
  which implies that
\begin{equation}
  \|J_2\|_{\mathscr{K}^{\sigma}_q}\lesssim\|f\|_{\mathscr{K}^{\sigma}_q}\|g\|_{\mathscr{K}^{\sigma}_q},
  \quad \forall f,g\in\mathscr{K}^{\sigma}_q(T), \quad \forall q\in [1,\infty],\;\; \forall \sigma\geqslant 0.
\end{equation}
  Combining (2.7) and (2.10), we obtain (2.3) in the case $1/2\leqslant\sigma\leqslant 1$ and $1/\sigma\leqslant q\leqslant 1/(1\!-\!\sigma)$.

  Proof of (2.3) in the rest case $\sigma>1$ and $1\leqslant q\leqslant\infty$ is much easier. Indeed, since for any $1<q\leqslant\infty$, the
  condition $\sigma>1$ implies that $2q'\sigma>1$, we see that for $\sigma>1$ and $1<q\leqslant\infty$,
\begin{eqnarray*}
  \sqrt{t}J_1(t)&\;\leqslant\;&\sqrt{2}\sup_{0<\tau<T}\sqrt{\tau}\Big|\ln\Big(\frac{\tau}{\rme T}\Big)\Big|^{\sigma}|f(\tau)|\cdot
  \int_0^t\Big|\ln\Big(\frac{\tau}{\rme T}\Big)\Big|^{-2\sigma}\cdot\sqrt{\tau}\Big|\ln\Big(\frac{\tau}{\rme T}\Big)\Big|^{\sigma}
  |g(\tau)|\frac{\rmd\tau}{\tau}
\\
  &\;\lesssim\;&\|f\|_{\dot{\mathscr{K}}^{\sigma}_{\infty}}\cdot\|g\|_{\dot{\mathscr{K}}^{\sigma}_q}
  \Big(\int_0^t\Big|\ln\Big(\frac{\tau}{\rme T}\Big)\Big|^{-2q'\sigma}\frac{\rmd\tau}{\tau}\Big)^{\frac{1}{q'}}
  \lesssim\|f\|_{\dot{\mathscr{K}}^{\sigma}_{\infty}}\|g\|_{\dot{\mathscr{K}}^{\sigma}_q}
  \Big|\ln\Big(\frac{t}{\rme T}\Big)\Big|^{1-\frac{1}{q}-2\sigma}.
\end{eqnarray*}
  Since $\sigma>1$ implies that $\displaystyle\Big|\ln\Big(\frac{t}{\rme T}\Big)\Big|^{1-\frac{1}{q}-\sigma}\in L^q\Big((0,T),\frac{dt}{t}\Big)$,
  we see that $J_1\in\dot{\mathscr{K}}^{\sigma}_q(T)$, and
\begin{equation}
  \|J_1\|_{\dot{\mathscr{K}}^{\sigma}_q}\leqslant\|f\|_{\dot{\mathscr{K}}^{\sigma}_{\infty}}\|g\|_{\dot{\mathscr{K}}^{\sigma}_q},
  \quad \forall f\in\mathscr{K}^{\sigma}_{\infty}(T),\;\; \forall g\in\mathscr{K}^{\sigma}_q(T).
\end{equation}
  It is easy to see that the above estimate also holds for the case $q=1$ and $\sigma>1$. From (2.10) and (2.11) we obtain (2.3) in the case $\sigma>1$
  and $1\leqslant q\leqslant\infty$. This completes the proof of Lemma 2.6. $\quad\Box$

\section{The proof of Theorem 1.2}
\setcounter{equation}{0}

\hskip 2em
  In this section we give the proof of Theorems 1.2. We shall first derive some linear and bilinear estimates, and next use these estimates to prove
  Theorem 1.2.

  Let $1\leqslant q<\infty$ and $\sigma\geqslant 0$. Given $T>0$, we introduce a path space $\mathscr{X}_T$ as follows:
\begin{eqnarray*}
  \mathscr{X}_T &=&\{\bfu\in L^{\infty}_{\rm loc}((0,T],L^\infty(\mathbb{R}^n)): \nabla\cdot\bfu=0,\;\;
  \|\bfu\|_{\mathscr{X}_T}<\infty\},
\end{eqnarray*}
  where
\begin{eqnarray*}
  \|\bfu\|_{\mathscr{X}_T}&=&\sup_{0<t<T}\sqrt{t}\Big|\ln\Big(\frac{t}{\rme T}\Big)\Big|^{\sigma}\|\bfu(t)\|_{\infty}
  +\Big[\int_0^T\!\!\Big(\sqrt{t}\Big|\ln\Big(\frac{t}{\rme T}\Big)\Big|^{\sigma}\|\bfu(t)\|_{\infty}\Big)^q\frac{\rmd t}{t}\,\Big]^{\frac{1}{q}}.
\end{eqnarray*}
  It is clear that $(\mathscr{X}_T,\|\cdot\|_{\mathscr{X}_T})$ is a Banach space. We shall also consider the following path spaces:
\begin{eqnarray*}
  &\mathscr{Y}_T=L^{\infty}((0,T),B^{-1,\sigma}_{\infty\, q}(\mathbb{R}^n))\cap\mathscr{X}_T,&
\\
  &\mathscr{Y}_T^0=C([0,T],B^{-1,\sigma}_{\infty\,q\,0}(\mathbb{R}^n))\cap\mathscr{X}_T.&
\end{eqnarray*}

  {\bf Lemma 3.1}\ \ {\em Let $1\leqslant q\leqslant\infty$ and $\sigma\geqslant 0$. If $\bfu_0\in B^{-1,\sigma}_{\infty\,q}(\mathbb{R}^n)$ then
  $e^{t\Delta}\bfu_0\in\mathscr{Y}_T$ for any finite $T>0$, and
$$
  \|e^{t\Delta}\bfu_0\|_{\mathscr{X}_T}+\sup_{t\in (0,T)}\|e^{t\Delta}\bfu_0\|_{B^{-1,\sigma}_{\infty\,q}}\lesssim_T
  \|\bfu_0\|_{B^{-1,\sigma}_{\infty\,q}}.
$$
  If furthermore $\bfu_0\in B^{-1,\sigma}_{\infty\,q\,0}(\mathbb{R}^n)$ then in addition to the above estimate we also have $e^{t\Delta}\bfu_0\in
  \mathscr{Y}_T^0$, and
$$
  \lim_{T\to 0^+}\|e^{t\Delta}\bfu_0\|_{\mathscr{X}_T}=0.
$$}

  {\em Proof}:\ \ It is easy to see that $B^{-1,\sigma}_{\infty\, q}(\mathbb{R}^n)$ is a shift-invariant Banach space of distributions. Hence by
  Propositions 4.1 and 4.4 of \cite{LEM02} we see that $\bfu_0\in B^{-1,\sigma}_{\infty\, q}(\mathbb{R}^n)$ implies that $e^{t\Delta}\bfu_0\in
  C_\ast([0,\infty),B^{-1,\sigma}_{\infty,\, q}(\mathbb{R}^n))$, i.e., for any $t\geqslant 0$ we have $e^{t\Delta}\bfu_0\in
  B^{-1,\sigma}_{\infty\, q}(\mathbb{R}^n)$, and the map $t\mapsto e^{t\Delta}\bfu_0$ from $[0,\infty)$ to $B^{-1,\sigma}_{\infty\, q}(\mathbb{R}^n)$
  is continuous for $t>0$ with respect to the norm topology of $B^{-1,\sigma}_{\infty\, q}(\mathbb{R}^n)$ and continuous at $t=0$ with respect to the
  $\ast$-weak topology of $B^{-1,\sigma}_{\infty\, q}(\mathbb{R}^n)$, and
$$
  \sup_{t>0}\|e^{t\Delta}\bfu_0\|_{B^{-1,\sigma}_{\infty\,q}}\leq\|\bfu_0\|_{B^{-1,\sigma}_{\infty\,q}}.
$$
  Moreover, from Lemma 2.1 (choosing $p=\infty$, $s=-1$, $\gamma=0$ and $t_0=T$) and the embedding $B^{-1,\sigma}_{\infty\, q}(\mathbb{R}^n)
  \subseteq B^{-1,\sigma}_{\infty\infty}(\mathbb{R}^n)$ (for $1\leqslant q<\infty$ and $\sigma\geqslant 0$) we see that $\bfu_0\in
  B^{-1,\sigma}_{\infty\, q}(\mathbb{R}^n)$ implies that $e^{t\Delta}\bfu_0\in\mathscr{X}_T$, and
$$
  \|e^{t\Delta}\bfu_0\|_{\mathscr{X}_T}\lesssim_T\|\bfu_0\|_{B^{-1,\sigma}_{\infty\,q}}.
$$
  Hence the first part of the lemma follows. The second part of the lemma follows from a standard density argument, cf. the proof of Lemma 2.5 in
  \cite{Cui}; we omit the details. $\quad\Box$
\medskip

  {\bf Lemma 3.2} \ \ {\em Let $T>0$ be given and assume that $(q,\sigma)$ satisfies one of the following two conditions: $(a)$ $\sigma\geqslant 1$
  and $1\leqslant q\leqslant\infty$; $(b)$ $1/2\leqslant\sigma<1$ and $1/\sigma\leqslant q\leqslant 1/(1\!-\!\sigma)$. Then $B(\bfu,\mathbf{v})\in
  {\mathscr{X}_T}$, and}
\begin{equation}
  \|B(\bfu,\mathbf{v})\|_{\mathscr{X}_T}\lesssim\|\bfu\|_{\mathscr{X}_T}\|\mathbf{v}\|_{\mathscr{X}_T}.
\end{equation}

  {\em Proof}:\ \ This is an immediate consequence of Lemma 2.2 and Lemma 2.5. $\quad\Box$
\medskip

  {\bf Lemma 3.3} \ \ {\em Let $T>0$ be given and assume that $(q,\sigma)$ satisfies one of the following two conditions: $(a)$ $\sigma\geqslant 1$
  and $1\leqslant q\leqslant\infty$; $(b)$ $1/2\leqslant\sigma<1$ and $1/\sigma\leqslant q\leqslant 1/(1\!-\!\sigma)$. Then for any $\bfu,\mathbf{v}
  \in {\mathscr{X}_T}$ we have $B(\bfu,\mathbf{v})\in L^{\infty}((0,T),B^{-1,\sigma}_{\infty\, q}(\mathbb{R}^n))$, and}
\begin{equation}
  \sup_{0<t<T}\|B(\bfu,\mathbf{v})(t)\|_{B^{-1,\sigma}_{\infty\, q}}\lesssim_T\|\bfu\|_{\mathscr{X}_T}\|\mathbf{v}\|_{\mathscr{X}_T}.
\end{equation}

  {\em Proof}:\ \ We first assume that $1\leqslant q<\infty$. By Lemma 2.2, for any $s>0$ and $t\in (0,T)$ we have
\begin{eqnarray*}
  \|e^{s\Delta}B(\bfu,\mathbf{v})(t)\|_{\infty}
  &=&\Big\|\int_0^t\! e^{(t+s-\tau)\Delta}\mathbb{P}\nabla\cdot[\bfu(\tau)\otimes\mathbf{v}(\tau)]\rmd\tau\Big\|_{\infty}
\\
  &\lesssim &\int_0^{t}\!(t+s-\tau)^{-\frac{1}{2}}\|\bfu(\tau)\|_{\infty}\|\mathbf{v}(\tau)\|_{\infty}\rmd\tau.
\end{eqnarray*}
  Hence
\begin{eqnarray*}
  && \Big[\int_0^T\!\!\Big(\sqrt{s}\Big|\ln\Big(\frac{s}{\rme T}\Big)\Big|^{\sigma}\|e^{s\Delta}B(\bfu,\mathbf{v})(t)\|_{\infty}
  \Big)^q\frac{\rmd s}{s}\Big]^{\frac{1}{q}}
\nonumber\\
  &\;\lesssim\; &\Big[\int_0^t\!\!\Big(\int_0^s\sqrt{s}\Big|\ln\Big(\frac{s}{\rme T}\Big)\Big|^{\sigma}(t+s-\tau)^{-\frac{1}{2}}
  \|\bfu(\tau)\|_{\infty}\|\mathbf{v}(\tau)\|_{\infty}\rmd\tau\Big)^q\frac{\rmd s}{s}\Big]^{\frac{1}{q}}
\nonumber\\
  && \quad +\Big[\int_0^t\!\!\Big(\int_s^t\sqrt{s}\Big|\ln\Big(\frac{s}{\rme T}\Big)\Big|^{\sigma}(t+s-\tau)^{-\frac{1}{2}}
  \|\bfu(\tau)\|_{\infty}\|\mathbf{v}(\tau)\|_{\infty}\rmd\tau\Big)^q\frac{\rmd s}{s}\Big]^{\frac{1}{q}}
\nonumber\\
  && \quad +\Big[\int_t^T\!\!\Big(\int_0^t\sqrt{s}\Big|\ln\Big(\frac{s}{\rme T}\Big)\Big|^{\sigma}(t+s-\tau)^{-\frac{1}{2}}
  \|\bfu(\tau)\|_{\infty}\|\mathbf{v}(\tau)\|_{\infty}\rmd\tau\Big)^q\frac{\rmd s}{s}\Big]^{\frac{1}{q}}
\nonumber\\ [0.3cm]
  &\;:=\; & K_1(t)+ K_2(t)+ K_3(t).
\end{eqnarray*}
  In $K_1(t)$ and $K_3(t)$ there holds the relation $\tau<s$, so that
$$
  K_1(t)+K_3(t)\lesssim\Big[\int_0^T\!\!\Big(\int_0^s\sqrt{s}\Big|\ln\Big(\frac{s}{\rme T}\Big)\Big|^{\sigma}(s-\tau)^{-\frac{1}{2}}
  \|\bfu(\tau)\|_{\infty}\|\mathbf{v}(\tau)\|_{\infty}\rmd\tau\Big)^q\frac{\rmd s}{s}\Big]^{\frac{1}{q}}, \quad \forall t\in (0,T).
$$
  By Lemma 2.5, the right-hand side is bounded by $\|\bfu\|_{\mathscr{X}_T}\|\mathbf{v}\|_{\mathscr{X}_T}$. Hence
\begin{equation}
  K_1(t)+K_3(t)\lesssim\|\bfu\|_{\mathscr{X}_T}\|\mathbf{v}\|_{\mathscr{X}_T}, \quad \forall t\in (0,T).
\end{equation}
  The estimate of $K_2(t)$ is easy. Indeed, by applying the Minkowsky inequality we have
\begin{eqnarray*}
  K_2(t)&\;\leqslant\;& \int_0^t\!\Big(\int_0^{\tau}s^{\frac{q}{2}-1}\Big|\ln\Big(\frac{s}{\rme T}\Big)\Big|^{\sigma q}\rmd s\Big)^{\frac{1}{q}}
  (t-\tau)^{-\frac{1}{2}}\|\bfu(\tau)\|_{\infty}\|\mathbf{v}(\tau)\|_{\infty}\rmd\tau
\nonumber\\
  &\lesssim &\int_0^t\!\sqrt{\tau}\Big|\ln\Big(\frac{\tau}{\rme T}\Big)\Big|^{\sigma}(t-\tau)^{-\frac{1}{2}}
  \|\bfu(\tau)\|_{\infty}\|\mathbf{v}(\tau)\|_{\infty}\rmd\tau
\nonumber\\
  &\lesssim &\sup_{0<\tau<T}\sqrt{\tau}\Big|\ln\Big(\frac{\tau}{\rme T}\Big)\Big|^{\sigma}\|\bfu(\tau)\|_{\infty}\cdot
  \sup_{0<\tau<T}\sqrt{\tau}\Big|\ln\Big(\frac{\tau}{\rme T}\Big)\Big|^{\sigma}\|\mathbf{v}(\tau)\|_{\infty}
\nonumber\\
  && \quad \times\int_0^t\!(t-\tau)^{-\frac{1}{2}}\tau^{-\frac{1}{2}}\Big|\ln\Big(\frac{\tau}{\rme T}\Big)\Big|^{-\sigma}\rmd\tau
  \lesssim\|\bfu\|_{\mathscr{X}_T}\|\mathbf{v}\|_{\mathscr{X}_T}, \quad \forall t\in (0,T).
\end{eqnarray*}
  Combining this estimate with (3.3) we see that
\begin{equation}
  \Big[\int_0^T\!\!\Big(\sqrt{s}\Big|\ln\Big(\frac{s}{\rme T}\Big)\Big|^{\sigma}\|e^{s\Delta}B(\bfu,\mathbf{v})(t)\|_{\infty}
  \Big)^q\frac{\rmd s}{s}\Big]^{\frac{1}{q}}\lesssim\|\bfu\|_{\mathscr{X}_T}\|\mathbf{v}\|_{\mathscr{X}_T}, \quad \forall t\in (0,T).
\end{equation}
  Similarly we can prove that
\begin{equation}
  \sup_{0<s<T}\sqrt{s}\Big|\ln\Big(\frac{s}{\rme T}\Big)\Big|^{\sigma}\|e^{s\Delta}B(\bfu,\mathbf{v})(t)\|_{\infty}
  \lesssim\|\bfu\|_{\mathscr{X}_T}\|\mathbf{v}\|_{\mathscr{X}_T}, \quad \forall t\in (0,T).
\end{equation}

  Having proved (3.4) and (3.5), we now apply Lemma 2.1 to conclude that $B(\bfu,\mathbf{v})(t)\in B^{-1,\sigma}_{\infty\, q}(\mathbb{R}^n)$ for all
  $0<t<T$, and moreover,
\begin{eqnarray*}
  \|B(\bfu,\mathbf{v})(t)\|_{B^{-1,\sigma}_{\infty\, q}}\lesssim_T\|\bfu\|_{\mathscr{X}_T}\|\mathbf{v}\|_{\mathscr{X}_T},
   \quad \forall t\in (0,T),
\end{eqnarray*}
  which proves (3.2). $\quad\Box$
\medskip

  {\bf Lemma 3.4} \ \ {\em Let $1\leqslant q<\infty$ and $\sigma\geqslant\sigma_q$. Let $T>0$ be given and assume that $\bfu,\mathbf{v}\in
  {\mathscr{X}_T}$. Then $B(\bfu,\mathbf{v})\in C_w([0,T],B^{-1,\sigma}_{\infty\, q}(\mathbb{R}^n))$, i.e., the map $t\mapsto B(\bfu,\mathbf{v})(t)$
  from $[0,T]$ to $B^{-1,\sigma}_{\infty\, q}(\mathbb{R}^n)$ is continuous with respect to $S'(\mathbb{R}^n)$-weak topology. If furthermore
  either $\bfu\in\mathscr{Y}_T^0$ or $\mathbf{v}\in\mathscr{Y}_T^0$ then also $B(\bfu,\mathbf{v})\in\mathscr{Y}_T^0$, and $B(\bfu,\mathbf{v})\in
  C([0,T],B^{-1,\sigma}_{\infty\, q0}(\mathbb{R}^n))$, i.e., $B(\bfu,\mathbf{v})(t)$ is continuous with respect to
  $B^{-1,\sigma}_{\infty\, q}(\mathbb{R}^n)$-norm, and moreover,}
\begin{equation}
  \lim_{t\to 0^+}\|B(\mathbf{u},\mathbf{v})(t)\|_{B^{-1,\sigma}_{\infty\,q}}=0.
\end{equation}

  {\em Proof}:\ \  We first prove that $B(\bfu,\mathbf{v})(t)$ is continuous at $t=0$ with respect to $S'(\mathbb{R}^n)$-weak topology. Indeed, since
  the condition $\sigma\geqslant\sigma_q$ implies that $\displaystyle\Big|\ln\Big(\frac{\tau}{\rme T}\Big)\Big|^{-2\sigma}\in L^{q'}\Big((0,T),
  \frac{\rmd\tau}{\tau}\Big)$, it follows that for any $\mathbf{w}\in S(\mathbb{R}^n)$ we have
\begin{eqnarray*}
  |\langle B(\bfu,\mathbf{v})(t),\mathbf{w}\rangle|
  &=&|\int_0^t\langle e^{(t-\tau)\Delta}\mathbb{P}\nabla\cdot[\bfu(\tau)\otimes\mathbf{v}(\tau)],\mathbf{w}\rangle \rmd\tau|
\\
  &\leqslant&\int_0^t\|\bfu(\tau)\|_{\infty}\|\mathbf{v}(\tau)\|_{\infty}\|\mathbb{P}\nabla\mathbf{w}\|_{1} \rmd\tau
\\
  &\leqslant&\|\bfu\|_{\mathscr{X}_T}\|\mathbb{P}\nabla\mathbf{w}\|_{1}\Big[\int_0^t\!\!\Big(\sqrt{\tau}
  \Big|\ln\Big(\frac{\tau}{\rme T}\Big)\Big|^{\sigma}\|\mathbf{v}(\tau)\|_{\infty}\Big)^q\frac{\rmd\tau}{\tau}\Big]^{\frac{1}{q}}
  \cdot\Big\|\Big|\ln\Big(\frac{\tau}{\rme T}\Big)\Big|^{-2\sigma}\Big\|_{L^{q'}((0,T),\frac{d\tau}{\tau})}
\\
  &\lesssim _T&\|\bfu\|_{\mathscr{X}_T}\|\mathbb{P}\nabla\mathbf{w}\|_{1}\Big[\int_0^t\!\!\Big(\sqrt{\tau}
  \Big|\ln\Big(\frac{\tau}{\rme T}\Big)\Big|^{\sigma}\|\mathbf{v}(\tau)\|_{\infty}\Big)^q\frac{\rmd\tau}{\tau}\Big]^{\frac{1}{q}}\to 0 \;\;
  (\mbox{as}\;\; t\to 0^+),
\end{eqnarray*}
  which proves the desired assertion. Next, let $0<t_0\leq T$. If $t_0<t<T$ then we write
\begin{eqnarray}
  &&B(\bfu,\mathbf{v})(t)-B(\bfu,\mathbf{v})(t_0)
\nonumber\\
  &=&\int_0^t e^{(t-\tau)\Delta}\mathbb{P}\nabla\cdot[\bfu(\tau)\otimes\mathbf{v}(\tau)]\rmd\tau
  -\int_0^{t_0}e^{(t_0-\tau)\Delta}\mathbb{P}\nabla\cdot[\bfu(\tau)\otimes\mathbf{v}(\tau)]\rmd\tau
\nonumber\\
  &=&\int_{t_0}^t e^{(t-\tau)\Delta}\mathbb{P}\nabla\cdot[\bfu(\tau)\otimes\mathbf{v}(\tau)]\rmd\tau
  +[e^{(t-t_0)\Delta}-I]\int_0^{t_0}e^{(t_0-\tau)\Delta}\mathbb{P}\nabla\cdot[\bfu(\tau)\otimes\mathbf{v}(\tau)]\rmd\tau
\nonumber\\
  &=:& A(t)+B(t),
\end{eqnarray}
  and if $0<t_0-\delta<t<t_0$ then we write
\begin{eqnarray}
  &&B(\bfu,\mathbf{v})(t_0)-B(\bfu,\mathbf{v})(t)
\nonumber\\
  &=&\int_0^{t_0}e^{(t_0-\tau)\Delta}\mathbb{P}\nabla\cdot[\bfu(\tau)\otimes\mathbf{v}(\tau)]\rmd\tau
  -\int_0^t e^{(t-\tau)\Delta}\mathbb{P}\nabla\cdot[\bfu(\tau)\otimes\mathbf{v}(\tau)]\rmd\tau
\nonumber\\
  &=&\int_t^{t_0} e^{(t_0-\tau)\Delta}\mathbb{P}\nabla\cdot[\bfu(\tau)\otimes\mathbf{v}(\tau)]\rmd\tau
  +[e^{(t_0-t)\Delta}-I]\int_0^t e^{(t-\tau)\Delta}\mathbb{P}\nabla\cdot[\bfu(\tau)\otimes\mathbf{v}(\tau)]\rmd\tau
\nonumber\\
  &=&\int_t^{t_0} e^{(t_0-\tau)\Delta}\mathbb{P}\nabla\cdot[\bfu(\tau)\otimes\mathbf{v}(\tau)]\rmd\tau
  +[e^{(t_0-t)\Delta}-I]\int_{t_0-\delta}^t e^{(t-\tau)\Delta}\mathbb{P}\nabla\cdot[\bfu(\tau)\otimes\mathbf{v}(\tau)]\rmd\tau
\nonumber\\
  &&+e^{(t-t_0+\delta)\Delta}[e^{(t_0-t)\Delta}-I]\int_0^{t_0-\delta}e^{(t_0-\delta-\tau)\Delta}\mathbb{P}\nabla\cdot
  [\bfu(\tau)\otimes\mathbf{v}(\tau)]\rmd\tau
\nonumber\\
  &=:& A_1(t)+B_1(t)+B_2(t).
\end{eqnarray}
  For $A(t)$ we have (see the proof of (2.7))
\begin{eqnarray*}
  \|A(t)\|_{B^{-1,\sigma}_{\infty\, q}}&\leqslant&\sup_{0<s<T}\sqrt{s}\Big|\ln\Big(\frac{s}{\rme T}\Big)\Big|^{\sigma}
  \Big\|\int_{t_0}^t e^{(t+s-\tau)\Delta}\mathbb{P}\nabla\cdot[\bfu(\tau)\otimes\mathbf{v}(\tau)]\rmd\tau\Big\|_{\infty}
\\
  &&+\Big[\int_0^T\!\!\Big(\sqrt{s}\Big|\ln\Big(\frac{s}{\rme T}\Big)\Big|^{\sigma}
  \Big\|\int_{t_0}^t e^{(t+s-\tau)\Delta}\mathbb{P}\nabla\cdot[\bfu(\tau)\otimes\mathbf{v}(\tau)]\rmd\tau\Big\|_{\infty}
  \Big)^q\frac{\rmd s}{s}\Big]^{\frac{1}{q}}
\\
  &\lesssim &
  \|\bfu\|_{\mathscr{X}_T}\|\mathbf{v}\|_{\mathscr{X}_T}\Big\{\sup_{0<s<T}\sqrt{s}\Big|\ln\Big(\frac{s}{\rme T}\Big)\Big|^{\sigma}
  \int_{t_0}^t(t+s-\tau)^{-\frac{1}{2}}\tau^{-1}\Big|\ln\Big(\frac{\tau}{\rme T}\Big)\Big|^{-2\sigma}\rmd\tau
\\
  &&+\Big[\int_0^T\!\!s^{\frac{q}{2}-1}\Big|\ln\Big(\frac{s}{\rme T}\Big)\Big|^{\sigma q}\Big(
  \Big\|\int_{t_0}^t(t+s-\tau)^{-\frac{1}{2}}\tau^{-1}\Big|\ln\Big(\frac{\tau}{\rme T}\Big)\Big|^{-2\sigma}\rmd\tau
  \Big)^q\rmd s\Big]^{\frac{1}{q}}\Big\}
\\
  &\lesssim_T&
  \|\bfu\|_{\mathscr{X}_T}\|\mathbf{v}\|_{\mathscr{X}_T}\int_{t_0}^t(t-\tau)^{-\frac{1}{2}}\tau^{-1}
  \Big|\ln\Big(\frac{\tau}{\rme T}\Big)\Big|^{-2\sigma}\rmd\tau.
\end{eqnarray*}
  Since $t_0>0$, we see that $\displaystyle\lim_{t\to t_0^+}\|A(t)\|_{B^{-1,\sigma}_{\infty\, q}}=0$. Moreover, since $B(t)=[e^{(t-t_0)\Delta}-I]
  B(\bfu,\mathbf{v})(t_0)$ and $B(\bfu,\mathbf{v})(t_0)\in B^{-1,\sigma}_{\infty\, q}(\mathbb{R}^n)$, by the assertion proved before we see that
  $\displaystyle\lim_{t\to t_0^+}B(t)=0$ in $S'(\mathbb{R}^n)$-weak topology. Hence
$$
  \lim_{t\to t_0^+}B(\bfu,\mathbf{v})(t)=B(\bfu,\mathbf{v})(t_0) \quad \mbox{in $S'(\mathbb{R}^n)$-weak topology}.
$$
  Next, similarly as for $A(t)$ we have $\displaystyle\lim_{t\to t_0^-}\|A_1(t)\|_{B^{-1,\sigma}_{\infty\, q}}=0$. Moreover, similarly as for the
  treatment of $A(t)$ we have that by choosing $\delta$ sufficiently small, $\|B_1(t)\|_{B^{-1,\sigma}_{\infty\, q}}$ can be as small as we expect,
  and when $\delta$ is chosen and fixed, $B_2(t)$ can be treated similarly as for $B(t)$ to get that for any $\mathbf{w}\in (S'(\mathbb{R}^n))^n$,
  $\displaystyle\lim_{t\to t_0^-}\langle B_2(t),\mathbf{w}\rangle=0$. Hence
$$
  \lim_{t\to t_0^-}B(\bfu,\mathbf{v})(t)=B(\bfu,\mathbf{v})(t_0) \quad \mbox{in $S'(\mathbb{R}^n)$-weak topology}.
$$
  This proves $B(\bfu,\mathbf{v})\in C_w([0,T],B^{-1,\sigma}_{\infty\, q}(\mathbb{R}^n))$. Finally, since $B(\bfu,\mathbf{v})\in
  L^{\infty}((0,T),B^{-1,\sigma}_{\infty\, q}(\mathbb{R}^n))$ for $\bfu,\bfv\in\mathscr{Y}_T$, it follows that if either $\bfu\in\mathscr{Y}_T^0$ or
  $\mathbf{v}\in\mathscr{Y}_T^0$ then by a standard density argument we see that $B(\bfu,\mathbf{v})\in C([0,T],B^{-1,\sigma}_{\infty\, q0}
  (\mathbb{R}^n))$, cf. the proof of the last assertion in Lemma 2.5 of \cite{Cui}. We omit the details here. $\quad\Box$
\medskip

  We are now ready to give the proof of Theorem 1.2.
\medskip

  {\bf Proof of Theorem 1.2}:\ \ Let $1\leqslant q<\infty$ and $\sigma\geqslant\sigma_q$ be given. We rewrite the problem (1.2) into the following
  equivalent integral equation:
$$
  \bfu(t)=e^{t\Delta}\bfu_0+B(\bfu,\bfu)(t).
$$
  Given $\bfu_0\in B^{-1,\sigma}_{\infty\, q}(\mathbb{R}^n)$ with $\diva\bfu_0=0$ and $T>0$, we define a map $\mathscr{J}:\mathscr{X}_T\to\mathscr{X}_T$ as
  follows: For any $\bfu\in\mathscr{X}_T$, $\mathscr{J}(\bfu)$ equals to the right-hand side of the above equation. By Lemmas 3.1 and 3.2, $\mathscr{J}$
  is a self-mapping in $\mathscr{X}_T$ and the following estimates hold:
\begin{eqnarray*}
  \|\mathscr{J}(\bfu)\|_{\mathscr{X}_T}&\,\leq\,&
  \|e^{t\Delta}\bfu_0\|_{\mathscr{X}_T}+C\|\bfu\|_{\mathscr{X}_T}^2,
\\
  \|\mathscr{J}(\bfu)-\mathscr{J}(\mathbf{v})\|_{\mathscr{X}_T}&\,\leq\,&
  C(\|\bfu\|_{\mathscr{X}_T}+\|\mathbf{v}\|_{\mathscr{X}_T})
  \|\bfu-\mathbf{v}\|_{\mathscr{X}_T}.
\end{eqnarray*}
  Choose a number $\varepsilon>0$ sufficiently small such that $4C\epsln<1$, where $C$ is the larger constant appearing in the above estimates. To prove
  the assertion $(1)$ of Theorem 1.2, for any $\bfu_0\in B^{-1,\sigma}_{\infty\,q\,0}(\mathbb{R}^n)$ with $\diva\bfu_0=0$ we choose $T>0$ so small that
  $\|e^{t\Delta}\bfu_0\|_{\mathscr{X}_T}\leq\epsln$. By Lemma 3.1, such $T$ exists. Then from Lemma 3.4 and the first inequality in the above we
  easily see that $\mathscr{J}$ maps the closed ball $\overline{B}(0,2\epsln)$ in $\mathscr{X}_T^0$ into itself, and the second inequality ensures that
  $\mathscr{J}$ is a contraction mapping when restricted to this ball. Hence, by the fixed point theorem of Banach, $\mathscr{J}$ has a unique fixed
  point in this ball. Since $e^{t\Delta}\bfu_0\in\mathscr{Y}_T^0$ and $B(\bfu,\bfu)\in\mathscr{Y}_T^0$ for $\bfu\in\mathscr{Y}_T^0$, from the iteration
  procedure we see that this fixed point lies in $\mathscr{Y}_T^0$. Hence we have obtain a mild solution of the problem (1.1) in the path space
  $\mathscr{Y}_T^0$. This proves the assertion $(1)$. To prove the assertion $(2)$, for given $T>0$ we let $\bfu_0\in B^{-1,\sigma}_{\infty\,q}
  (\mathbb{R}^n)$ (with $\diva\bfu_0=0$) be so small that $\|e^{t\Delta}\bfu_0\|_{\mathscr{X}_T}\leq\epsln$. Then from the first inequality in the above
  we easily see that $\mathscr{J}$ maps the closed ball $\overline{B}(0,2\epsln)$ in $\mathscr{X}_T$ into itself, and the second inequality ensures that
  $\mathscr{J}$ is a contraction mapping when restricted to this ball. Hence, again by the fixed point theorem of Banach, $\mathscr{J}$ has a unique
  fixed point in this ball. Since $e^{t\Delta}\bfu_0\in\mathscr{Y}_T$ and $B(\bfu,\bfu)\in\mathscr{Y}_T$ for $\bfu\in \mathscr{Y}_T$, by a similar
  argument as above we get a mild solution of the problem (1.1) which lies in the path space $\mathscr{Y}_T$. This proves the assertion $(2)$. The
  proof of Theorem 1.2 is complete. $\quad\Box$

\section{The proof of Theorem 1.3}
\setcounter{equation}{0}

\hskip 2em
  In this section we give the proof of Theorem 1.3. We shall mainly consider the case $2<q\leqslant\infty$ and $1-2/q\leqslant\sigma<1-1/q$, because
  the rest cases $1\leqslant q\leqslant 2$, $0\leqslant\sigma<1/q$ and $2<q\leqslant \infty$, $0\leqslant\sigma<1-2/q$ are easier to treat. In the
  end of this section we shall explain how to modify the arguments given below to get proofs for these two cases. Besides, we only give proof for the
  case $n\geqslant 3$, and proof for the two-dimension case is omitted.

  For a sufficiently large positive integer $m$, we denote
$$
  \mathscr{A}_m=\{4k:\;k\in\mathbb{N},\; 4m+1\leqslant k\leqslant 5m\}, \qquad
  \mathscr{B}_m=\{4k:\;k\in\mathbb{N},\; m+1\leqslant k\leqslant 2m\}.
$$
  Clearly $|\mathscr{A}_m|=|\mathscr{B}_m|=m$. Let $\varepsilon$ be a sufficiently small positive number. For every positive integer $k$, we introduce
  three $n$-dimensional vector $a_k$, $b_k$ and $c_k$ as follows:
$$
  a_k=2^kn^{-\frac{1}{2}}(1,1,\cdots,1), \quad b_k=2^{k-1}(\varepsilon,2\varepsilon,\sqrt{1-5\varepsilon^2},0,\cdots,0), \quad
  c_k=2^k(1,0,\cdots,0).
$$
  Note that $|a_k|=2|b_k|=|c_k|=2^k$, $k=1,2,\cdots$. Let $\phi$ be as in Section 2 and set $\rho(\xi)=\phi(8|\xi|)$, $\xi\in\mathbb{R}^n$. It is clear
  that $\rho(\xi)=1$ for $|\xi|\leqslant 5/32$ and $\supp\rho\subseteq\bar{B}(0,3/16)$. We denote
\begin{eqnarray*}
  \Phi_{kl}^{++}(\xi)=\rme^{\bfi c_l\xi}\rho(\xi-a_k-b_l),  & \quad &   \Phi_{kl}^{+-}(\xi)=\rme^{\bfi c_l\xi}\rho(\xi-a_k+b_l),
\\
  \Phi_{kl}^{-+}(\xi)=\rme^{\bfi c_l\xi}\rho(\xi+a_k-b_l),  & \quad &   \Phi_{kl}^{--}(\xi)=\rme^{\bfi c_l\xi}\rho(\xi+a_k+b_l),
\end{eqnarray*}
$$
  \Psi_{kl}(\xi)=\Phi_{kl}^{++}(\xi)+\Phi_{kl}^{+-}(\xi)+\Phi_{kl}^{-+}(\xi)+\Phi_{kl}^{--}(\xi)
$$
  ($k,l=1,2,\cdots$). We now consider the initial value problem
\begin{eqnarray}
\left\{
\begin{array}{l}
  \partial_t\bfu-\Delta\bfu+\mathbb{P}\nabla\cdot(\bfu\otimes\bfu)=0\quad \mbox{in}\;\,\mathbb{R}^n\times\mathbb{R}_+,\\
  \bfu(x,t)=\delta\bfu_0(x)\quad \mbox{for}\;\,  x\in\mathbb{R}^n,
\end{array}
\right.
\end{eqnarray}
  with a sufficiently small $\delta>0$ and $\bfu_0=(u_1^0,u_2^0,\cdots,u_n^0)$, where
\begin{equation}
\left\{
\begin{array}{rcl}
  u_1^0(x)&=&\displaystyle m^{-\sigma-\frac{1}{q}}\sum_{k\in\mathscr{A}_m}\sum_{l\in\mathscr{B}_m}2^k\mathscr{F}^{-1}[\Psi_{kl}(\xi)],
\\ [0.3cm]
  u_2^0(x)&=&\displaystyle -m^{-\sigma-\frac{1}{q}}\sum_{k\in\mathscr{A}_m}\sum_{l\in\mathscr{B}_m}
  2^k\mathscr{F}^{-1}\Big[\frac{\xi_1}{\xi_2}\Psi_{kl}(\xi)\Big],
\\ [0.3cm]
  u_3^0(x)&=&\cdots=u_n^0(x)=0
\end{array}
\right.
\end{equation}
  (for $x\in\mathbb{R}^n$). here we follow the convention that $1/\infty=0$. Note that
$$
  u_1^0(x)= m^{-\sigma-\frac{1}{q}}\sum_{k\in\mathscr{A}_m}\sum_{l\in\mathscr{B}_m}2^k[\cos(a_k+b_l)(x+c_l)+\cos(a_k-b_l)(x+c_l)]
  \check{\rho}(x+c_l)
$$
  and
$$
  u_2^0(x)=-u_1^0(x)+m^{-\sigma-\frac{1}{q}}\sum_{k\in\mathscr{A}_m}\sum_{l\in\mathscr{B}_m}
  2^k\mathscr{F}^{-1}\Big[\frac{\xi_2-\xi_1}{\xi_2}\Psi_{kl}(\xi)\Big].
$$
 Note also that $u_i^0\in S'(\mathbb{R}^n)$ ($i=1,2,\cdots,n$) and $\diva\bfu_0=0$. Besides, it is clear that
\begin{equation}
  \supp\Psi_{kl}\subseteq \{\xi\in\mathbb{R}^n:2^{k-1}<|\xi|<2^{k+1}\} \quad \mbox{for}\;\; k\in\mathscr{A}_m, \;\; l\in\mathscr{B}_m.
\end{equation}

  {\bf Lemma 4.1} \ \ {\em Let $m\gg-\ln\varepsilon$. Then for any $1\leqslant q\leqslant\infty$ we have}
\begin{equation}
  \|u_i^0\|_{B^{-1,\sigma}_{\infty\,q}}\lesssim 1, \quad i=1,2.
\end{equation}

  {\em Proof}:\ \ By (4.1), it is clear that $\supp\,\widehat{u_1^0}$ does not intersects $\bar{B}(0,3/2)$. Hence
\begin{eqnarray*}
  \|u_1^0\|_{B^{-1,\sigma}_{\infty\,q}}&=&\Big[\sum_{j=1}^{\infty}\Big(2^{-j}j^{\sigma}\|\check{\psi}_j\ast u_1^0\|_{\infty}\Big)^q\,\Big]^{\frac{1}{q}}
  \lesssim m^{-\sigma-\frac{1}{q}}\Big[\sum_{j\in\mathscr{A}_m}
  \Big(j^{\sigma}\Big\|\sum_{l\in\mathscr{B}_m}|\check{\rho}(x+c_l)|\Big\|_{\infty}\Big)^q\,\Big]^{\frac{1}{q}}
\\
  &\lesssim\; & m^{-\sigma-\frac{1}{q}}\cdot m^{\sigma+\frac{1}{q}}\Big\|\sum_{l\in\mathscr{B}_m}|\check{\rho}(x+c_l)|\Big\|_{\infty}
  \lesssim 1 \quad \mbox{for}\;\; m\gg-\ln\varepsilon.
\end{eqnarray*}
  Next, choose a function $\chi\in C^{\infty}(\mathbb{R}^n\backslash\{0\})$ such that it is homogeneous of degree zero, $\chi(\xi)=1$ for
  $|\xi/|\xi|-e|\leqslant 1/32$, where $e=(1,1,\cdots,1)/\sqrt{n}$, and $\chi(\xi)=0$ for $|\xi/|\xi|-e|\geqslant 1/16$, and set $\psi'(\xi)=
  (\xi_1/\xi_2)\chi(\xi)\psi(\xi)$, $\psi'_j(\xi)=\psi'(2^{-j}\xi)$, $j=1,2,\cdots$. Then since $\psi_j(\xi)\widehat{u_2^0}(\xi)=\psi'_j(\xi)
  \widehat{u_1^0}(\xi)$, $j=1,2,\cdots$, we have
\begin{eqnarray*}
  \|u_2^0\|_{B^{-1,\sigma}_{\infty\,q}}&=&\Big[\sum_{j=1}^{\infty}\Big(2^{-j}j^{\sigma}\|\check{\psi}_j\ast u_2^0\|_{\infty}\Big)^q\,\Big]^{\frac{1}{q}}
  =\Big[\sum_{j=1}^{\infty}\Big(2^{-j}j^{\sigma}\|\check{\psi}_j'\ast u_1^0\|_{\infty}\Big)^q\,\Big]^{\frac{1}{q}}
  \lesssim 1.
\end{eqnarray*}
  The last inequality follows from a similar argument as above. This proves the lemma. $\quad\Box$
\medskip

  In what follows, for $s\in\mathbb{R}$, $\sigma\geqslant 0$, $p,q\in[1,\infty]$ and nonempty subset $A$ of $\mathbb{N}$, we denote
$$
  \|u\|_{B^{s,\sigma}_{p\,q}(A)}=\left\{
\begin{array}{ll}
  \displaystyle\Big[\sum_{j\in A}\Big(2^{js}j^{\sigma}\|\Delta_j u\|_p\Big)^q\Big]^{\frac{1}{q}} \quad &\mbox{for}\;\; 1\leqslant q<\infty,
\\ [0.3cm]
  \displaystyle\sup_{j\in A}\Big(2^{js}j^{\sigma}\|\Delta_j u\|_p\Big) \quad &\mbox{for}\;\; q=\infty.
\end{array}
\right.
$$

  {\bf Lemma 4.2} \ \ {\em Let $m\gg-\ln\varepsilon$ and $t=\varepsilon2^{-32m}$. Then for any $1\leqslant q\leqslant\infty$ we have}
\begin{equation}
  \Big\|\int_0^t\rme^{(t-\tau)\Delta}(\partial_1-\partial_2)(\rme^{\tau\Delta}u_1^0\rme^{\tau\Delta}u_1^0)\rmd\tau
  \Big\|_{B^{-1,\sigma}_{\infty\,q}(\mathscr{B}_m)}\gtrsim\varepsilon^2 m^{1-\sigma-\frac{1}{q}}.
\end{equation}

  {\em Proof}:\ \ First we note that for any $f\in S'(\mathbb{R}^n)$, if we denote $u=\rme^{t\Delta}f$ then $\supp\hat{u}(\cdot,t)=
  \supp\hat{f}$ for all $t>0$. Hence, if the frequency support of $f$ satisfies certain property, then the same property is also
  satisfied by $u=\rme^{t\Delta}f$ for all $t>0$.

  We shall use the following principle to prove the above result: If $f\in C^1(\mathbb{R}^n)$ is a real valued function and
  $\nabla f\in L^{\infty}(\mathbb{R}^n)$, then for any $\nu,x_0,x_0'\in\mathbb{R}^n$ with $|\nu|=1$, $x_0\neq x_0'$ and
  $x_0-x_0'\parallel\nu$, there holds
$$
  \|\partial_{\nu}f\|_{L^{\infty}(\mathbb{R}^n)}\geqslant\frac{|f(x_0)-f(x_0')|}{|x_0-x_0'|}.
$$
  In what follows, we shall use this principle to the case $\nu=(1/\sqrt{2})(1,-1,0,\cdots,0)$ (so that $\partial_{\nu}
  =(1/\sqrt{2})(\partial_1-\partial_2)$) and
\begin{equation}
  f(x)=\int_0^t\rme^{(t-\tau)\Delta}U_{1j}\rmd\tau
\end{equation}
  (see (4.7) for the expression of the function $U_{1j}$).

  Since $\supp\rho(\cdot-a)\ast\rho(\cdot-b)\subseteq B(a+b,1/2)$, we see that for any $k,k'\in\mathscr{A}_m$ and
  $l,l'\in\mathscr{B}_m$,
$$
  \supp(\Phi_{kl}^{+\mu}\ast\Phi_{kl'}^{+\nu})\subseteq B(2a_k,2^{\frac{k}{2}+1}), \quad \mu,\nu\in\{0,1\},
$$
$$
  \supp(\Phi_{kl}^{-\mu}\ast\Phi_{kl'}^{-\nu})\subseteq B(-2a_k,2^{\frac{k}{2}+1}), \quad \mu,\nu\in\{0,1\},
$$
$$
  \supp(\Phi_{kl}^{++}\ast\Phi_{kl}^{--}),\;\supp(\Phi_{kl}^{+-}\ast\Phi_{kl}^{-+})\subseteq B(0,1),
$$
$$
  \supp(\Phi_{kl}^{++}\ast\Phi_{kl'}^{--}),\;\supp(\Phi_{kl}^{+-}\ast\Phi_{kl'}^{-+})\subseteq
  \{\xi\in\mathbb{R}^n:\;2^{l\wedge l'-2}\leqslant|\xi|<2^{l\wedge l'-1}\} \;\;  (l\neq l'),
$$
$$
  \supp(\Phi_{kl}^{\mu\nu}\ast\Phi_{k'l'}^{\mu'\nu'})\subseteq \{\xi\in\mathbb{R}^n:\;2^{k\wedge k'-1}<|\xi|<2^{k\wedge k'+1}\}
  \;\;  (k\neq k'), \quad \mu,\nu\in\{0,1\},
$$
  where as usual $k\wedge k'=\max\{k,k'\}$. Moreover, it is easy to see that for any $j,l,l'\in\mathscr{B}_m$, if $l\wedge l'\neq j$ then
$$
  \Delta_j[\mathscr{F}^{-1}(\Phi_{kl}^{++}\ast\Phi_{kl'}^{-+})]=\Delta_j[\mathscr{F}^{-1}(\Phi_{kl}^{+-}\ast\Phi_{kl'}^{--})]=0,
$$
  and if $l\wedge l'=j$ then
$$
  \Delta_j[\mathscr{F}^{-1}(\Phi_{kl}^{++}\ast\Phi_{kl'}^{-+})]=
\left\{
\begin{array}{l}
  \mathscr{F}^{-1}(\Phi_{kj}^{++}\ast\Phi_{kj'}^{-+}) \quad \mbox{if}\;\, j=l,\\
  e^{2{\bfi}a_k(c_{j'}-c_j)}\mathscr{F}^{-1}(\Phi_{kj}^{++}\ast\Phi_{kj'}^{-+}) \quad \mbox{if}\;\, j=l',
\end{array}
\right.
$$
$$
  \Delta_j[\mathscr{F}^{-1}(\Phi_{kl}^{+-}\ast\Phi_{kl'}^{--})]=
\left\{
\begin{array}{l}
  \mathscr{F}^{-1}(\Phi_{kj}^{+-}\ast\Phi_{kj'}^{--}) \quad \mbox{if}\;\, j=l,\\
  e^{2{\bfi}a_k(c_{j'}-c_j)}\mathscr{F}^{-1}(\Phi_{kj}^{+-}\ast\Phi_{kj'}^{--}) \quad \mbox{if}\;\, j=l',
\end{array}
\right.
$$
  where $j'=l\vee l'=\min\{l,l'\}$. It follows that for any $j\in\mathscr{B}_m$,
\begin{eqnarray}
  &&\Delta_j(\rme^{\tau\Delta}u_1^0\rme^{\tau\Delta}u_1^0)
\nonumber\\
  &= & m^{-2\sigma-\frac{2}{q}}\sum_{k\in\mathscr{A}_m}\sum_{l,l'\in\mathscr{B}_m}2^{2k+1}
  \Delta_j[\mathscr{F}^{-1}(\rme^{-\tau|\xi|^2}\Phi_{kl}^{++}\ast\rme^{-\tau|\xi|^2}\Phi_{kl'}^{-+}
  +\rme^{-\tau|\xi|^2}\Phi_{kl}^{+-}\ast\rme^{-\tau|\xi|^2}\Phi_{kl'}^{--})]
\nonumber\\
  &= & m^{-2\sigma-\frac{2}{q}}\sum_{k\in\mathscr{A}_m}2^{2k+1}\mathscr{F}^{-1}
  (\rme^{-\tau|\xi|^2}\Phi_{kj}^{++}\ast\rme^{-\tau|\xi|^2}\Phi_{kj}^{-+}
  +\rme^{-\tau|\xi|^2}\Phi_{kj}^{+-}\ast\rme^{-\tau|\xi|^2}\Phi_{kj}^{--})
\nonumber\\
  && +m^{-2\sigma-\frac{2}{q}}\sum_{k\in\mathscr{A}_m}\sum_{{j'\in\mathscr{B}_m\atop j'<j}}\theta_{kjj'}2^{2k+1}
  \Delta_j[\mathscr{F}^{-1}(\rme^{-\tau|\xi|^2}\Phi_{kj}^{++}\ast\rme^{-\tau|\xi|^2}\Phi_{kj'}^{-+}
  +\rme^{-\tau|\xi|^2}\Phi_{kj}^{+-}\ast\rme^{-\tau|\xi|^2}\Phi_{kj'}^{--})]
\nonumber\\
  &\;:=\; & m^{-2\sigma-\frac{2}{q}}U_{1j}+m^{-2\sigma-\frac{2}{q}}U_{2j},
\end{eqnarray}
  where $\theta_{kjj'}=1+e^{2{\bfi}a_k(c_{j'}-c_j)}$. Hence
\begin{eqnarray}
  &&\Big\|\int_0^t\rme^{(t-\tau)\Delta}(\partial_1-\partial_2)(\rme^{\tau\Delta}u_1^0\rme^{\tau\Delta}u_1^0)\rmd\tau
  \Big\|_{B^{-1,\sigma}_{\infty\,q}(\mathscr{B}_m)}
\nonumber\\
  &\gtrsim\; &m^{-2\sigma-\frac{2}{q}}\Big[\sum_{j\in\mathscr{B}_m}\Big(2^{-j}j^{\sigma}
  \Big\|(\partial_1-\partial_2)\int_0^t\rme^{(t-\tau)\Delta}U_{1j}\rmd\tau\Big\|_{\infty}\Big)^q\;\Big]^{\frac{1}{q}}
\nonumber\\
  &&\quad -m^{-2\sigma-\frac{2}{q}}\Big[\sum_{j\in\mathscr{B}_m}\Big(2^{-j}j^{\sigma}
  \Big\|(\partial_1-\partial_2)\int_0^t\rme^{(t-\tau)\Delta}U_{2j}\rmd\tau\Big\|_{\infty}\Big)^q\;\Big]^{\frac{1}{q}}
\nonumber\\ [0.2cm]
  &:=\; &I+I\!\!I.
\end{eqnarray}
  Let $f$ be the function as defined in (4.6). We have
\begin{eqnarray*}
  f(x)&\;=\;&\int_0^t\rme^{(t-\tau)\Delta}U_{1j}\rmd\tau
\\
  &\;=\;&\mathscr{F}^{-1}\Big[\sum_{k\in\mathscr{A}_m}2^{2k+1}\int_0^t\rme^{-(t-\tau)|\xi|^2}\Big(\rme^{-\tau|\xi|^2}\Phi_{kj}^{++}
  \ast\rme^{-\tau|\xi|^2}\Phi_{kj}^{-+}+\rme^{-\tau|\xi|^2}\Phi_{kj}^{+-}\ast\rme^{-\tau|\xi|^2}\Phi_{kj}^{--}\Big)\rmd\tau\Big]
\\
  &\;=\;&\mathscr{F}^{-1}\Big[\sum_{k\in\mathscr{A}_m}2^{2k+1}\rme^{\bfi c_j\xi}\!\int_0^t\rme^{-(t-\tau)|\xi|^2}
  \int_{\mathbb{R}^n}\rme^{-\tau(|\xi-\eta|^2+|\eta|^2)}
\\
  &&\quad \times\Big(\rho(\xi-\eta-a_k-b_j)\rho(\eta+a_k-b_j)+\rho(\xi-\eta-a_k+b_j)\rho(\eta+a_k+b_j)\Big)\rmd\eta\rmd\tau\Big]
\\
   &\;=\;&\mathscr{F}^{-1}\{\rme^{\bfi c_j\xi}[g(t,\xi)+g(t,-\xi)]\}
\\
  &\;=\;&\frac{2}{(2\pi)^n}\int_{\mathbb{R}^n}\cos[(x+c_j)\xi]g(t,\xi)\rmd\xi,
\end{eqnarray*}
  where
$$
  g(t,\xi)=\rme^{-t|\xi|^2}\int_{\mathbb{R}^n}G_{j}(t,\xi,\eta)\rho(\xi-\eta-2b_j)\rho(\eta)\rmd\eta,
$$
  and
$$
  G_{j}(t,\xi,\eta)=\sum_{k\in\mathscr{A}_m}2^{2k+1}\int_0^t\rme^{-\tau(|\xi-\eta+a_k-b_j|^2+|\eta-a_k+b_j|^2-|\xi|^2)}\rmd\tau.
$$
  Let $x_0=-c_j$, $x_0'=-c_j+2^{-j-1}\pi\varepsilon^{-1}(1,-1,0,\cdots,0)$. Then since $x_0-x_0'\parallel (1,-1,0,\cdots,0)$ and
  $|x_0-x_0'|=2^{-j-\frac{1}{2}}\pi\varepsilon^{-1}\sim 2^{-j}\varepsilon^{-1}$, we have
\begin{eqnarray*}
  \Big\|(\partial_1-\partial_2)\int_0^t\rme^{(t-\tau)\Delta}U_{1j}\rmd\tau\Big\|_{\infty}
  &\;\gtrsim\;&\frac{1}{2^{-j}\varepsilon^{-1}}\Big|\int_{\mathbb{R}^n}\{1-\cos[2^{-j-1}\pi\varepsilon^{-1}(\xi_1-\xi_2)]\}
  g(t,\xi)\rmd\xi\Big|.
\end{eqnarray*}
  From the expression of $g(t,\xi)$ we see that on the support of $g$ we have $\xi\sim 2b_j$, which implies that
  $2^{-j-1}\pi\varepsilon^{-1}(\xi_1-\xi_2)\sim \pi/2$, and, consequently,
$$
  1-\cos[2^{-j-1}\pi\varepsilon^{-1}(\xi_1-\xi_2)]\sim 1.
$$
  Hence
\begin{eqnarray*}
  \Big\|(\partial_1-\partial_2)\int_0^t\rme^{(t-\tau)\Delta}U_{1j}\rmd\tau\Big\|_{\infty}
  &\;\gtrsim\;&\frac{1}{2^{-j}\varepsilon^{-1}}\int_{\mathbb{R}^n}g(t,\xi)\rmd\xi.
\end{eqnarray*}
  Since $\supp\rho\subseteq\bar{B}(0,3/4)$, and
$$
  |a_k|=2^k\geqslant 2^{16m+4}, \quad  |b_j|=2^{j-1}\leqslant 2^{8m}, \quad \forall k\in\mathscr{A}_m,\;\; \forall j\in\mathscr{B}_m,
$$
  we see that on the support of $\rho(\xi-\eta-2b_j)\rho(\eta)$ there holds
$$
  |\xi-\eta+a_k-b_j|^2+|\eta-a_k+b_j|^2-|\xi|^2\sim 2|a_k|^2=2^{2k+1},
$$
  so that
$$
  G_{j}(t,\xi,\eta)\sim\sum_{k\in\mathscr{A}_m}2^{2k+1}\frac{1-\rme^{-t2^{2k+1}}}{2^{2k+1}}\gtrsim m(1-\rme^{-t2^{32m}}).
$$
  Hence,
\begin{eqnarray*}
  \int_{\mathbb{R}^n}g(t,\xi)\rmd\xi &\;\gtrsim\;&\rme^{-t2^{2j}}\cdot m(1-\rme^{-t2^{32m}}),
\end{eqnarray*}
  and, consequently,
\begin{eqnarray*}
  \Big\|(\partial_1-\partial_2)\int_0^t\rme^{(t-\tau)\Delta}U_{1j}\rmd\tau\Big\|_{\infty}
  &\;\gtrsim\;&2^j\varepsilon\rme^{-t2^{2j}}\cdot m(1-\rme^{-t2^{32m}}).
\end{eqnarray*}
  In getting the last inequality we have used the assumption that $t=\varepsilon 2^{-32m}$ and $0<\varepsilon\ll 1$. It follows that
\begin{equation}
  I\gtrsim m^{-2\sigma-\frac{2}{q}}\cdot\varepsilon^2 m^{1+\sigma+\frac{1}{q}}\gtrsim\varepsilon^2 m^{1-\sigma-\frac{1}{q}}.
\end{equation}
  Next, for any $f\in L^1(\mathbb{R}^n)$, by writing
$$
  \rme^{t\Delta}f(x)=\int_{\mathbb{R}^n}\rme^{-|y|^2}f(x-2\sqrt{t}y)\rmd y, \quad x\in\mathbb{R}^n, \;\; t>0,
$$
  we see that if $f$ satisfies the property
$$
  |f(x)|\leqslant C_N(1+|x|)^{-N}, \quad \forall x\in\mathbb{R}^n, \;\; \forall N>0,
$$
  then there also holds
$$
  |\rme^{t\Delta}f(x+a)|\leqslant C_N(1+|x+a|)^{-N}, \quad \forall x\in\mathbb{R}^n, \;\; \forall t\in (0,1),
   \;\; \forall N>0, \;\; \forall a\in\mathbb{R}^n.
$$
  It follows that for any $k\in\mathscr{A}_m$ and $j,j'\in\mathscr{B}_m$ with $j'<j$ we have
\begin{eqnarray*}
  \|\partial_i\Delta_j(\rme^{\tau\Delta}\check{\Phi}_{kj}^{+\mu}\cdot\rme^{\tau\Delta}\check{\Phi}_{kj'}^{-\mu})\|_{\infty}
  &\;\lesssim\;& 2^j\|\rme^{\tau\Delta}\check{\Phi}_{kj}^{+\mu}\cdot\rme^{\tau\Delta}\check{\Phi}_{kj'}^{-\mu}\|_{\infty}
\nonumber\\
  &\;\lesssim_N\;& 2^j\|(1+|x+c_j|)^{-N}(1+|x+c_{j'}|)^{-N}\|_{\infty}
\nonumber\\
  &\;\lesssim_N\;& 2^{-(N-1)j}, \quad  \forall N>0, \;\; \mu\in\{+,-\}, \;\; i=1,2.
\end{eqnarray*}
  The last estimate follows from a similar argument as in the proof of (2.45) of \cite{Wang}. Hence
\begin{eqnarray*}
  \Big\|\int_0^t\rme^{(t-\tau)\Delta}(\partial_1-\partial_2)U_{2j}\rmd\tau\Big\|_{\infty}
  &\;\lesssim\;& \Big(\sum_{k\in\mathscr{A}_m}2^{2k+1}\Big)\cdot tm2^{-(N-1)j}
\\
  &\lesssim & 2^{40m}\cdot tm2^{-4m(N-1)}\lesssim\varepsilon m,
\end{eqnarray*}
  where we have put $N=3$. It follows that
\begin{equation}
  I\!\!I\lesssim m^{-2\sigma-\frac{2}{q}}\cdot 2^{-4m}m^{\sigma}\cdot\varepsilon m\lesssim\varepsilon m^{1-\sigma-\frac{2}{q}}2^{-4m}.
\end{equation}
  Substituting (4.9) and (4.10) into (4.8), and assuming that $m$ is so large that $2^{-4m}\leqslant\varepsilon^2$, we obtain (4.5). $\quad\Box$
\medskip

  {\bf Lemma 4.3} \ \ {\em Let $m\gg-\ln\varepsilon$. Then for any $1\leqslant q\leqslant\infty$ and $t>0$ we have}
\begin{equation}
  \Big\|\int_0^t\rme^{(t-\tau)\Delta}\partial_2[\rme^{\tau\Delta}(u_1^0+u_2^0)\rme^{\tau\Delta}u_1^0]\rmd\tau
  \Big\|_{B^{-1,\sigma}_{\infty\,q}(\mathscr{B}_m)}\lesssim \varepsilon^2 2^{-8m}m^{-\sigma-\frac{2}{q}}.
\end{equation}

  {\em Proof}:\ \ Similarly as in the proof of the above lemma, for any $j\in\mathscr{B}_m$ we have
\begin{eqnarray*}
  &&\Delta_j[\rme^{\tau\Delta}(u_1^0+u_2^0)\rme^{\tau\Delta}u_1^0]
\\
  &= & m^{-2\sigma-\frac{2}{q}}\sum_{k\in\mathscr{A}_m}\sum_{l,l'\in\mathscr{B}_m}2^{2k+1}
  \Delta_j[\mathscr{F}^{-1}(\rme^{-\tau|\xi|^2}\frac{\xi_2-\xi_1}{\xi_2}\Phi_{kl}^{++}\ast\rme^{-\tau|\xi|^2}\Phi_{kl'}^{-+}
\\
  && \qquad\qquad\qquad\qquad\qquad\qquad\quad +\rme^{-\tau|\xi|^2}\frac{\xi_2-\xi_1}{\xi_2}\Phi_{kl}^{+-}\ast\rme^{-\tau|\xi|^2}\Phi_{kl'}^{--})]
\\
  &= & m^{-2\sigma-\frac{2}{q}}\sum_{k\in\mathscr{A}_m}2^{2k+1}\mathscr{F}^{-1}
  (\rme^{-\tau|\xi|^2}\frac{\xi_2-\xi_1}{\xi_2}\Phi_{kj}^{++}\ast\rme^{-\tau|\xi|^2}\Phi_{kj}^{-+}
  +\rme^{-\tau|\xi|^2}\frac{\xi_2-\xi_1}{\xi_2}\Phi_{kj}^{+-}\ast\rme^{-\tau|\xi|^2}\Phi_{kj}^{--})
\\
  && +m^{-2\sigma-\frac{2}{q}}\sum_{k\in\mathscr{A}_m}\sum_{{j'\in\mathscr{B}_m\atop j'<j}}\theta_{kjj'}2^{2k+1}
  \mathscr{F}^{-1}(\rme^{-\tau|\xi|^2}\frac{\xi_2-\xi_1}{\xi_2}\Phi_{kj}^{++}\ast\rme^{-\tau|\xi|^2}\Phi_{kj'}^{-+}
\\
  && \qquad\qquad\qquad\qquad\qquad\qquad\qquad +\rme^{-\tau|\xi|^2}\frac{\xi_2-\xi_1}{\xi_2}\Phi_{kj}^{+-}\ast\rme^{-\tau|\xi|^2}\Phi_{kj'}^{--})
\\
  &\;:=\; & m^{-2\sigma-\frac{2}{q}}V_{1j}+m^{-2\sigma-\frac{2}{q}}V_{2j}.
\end{eqnarray*}
  Hence
\begin{eqnarray}
  &&\Big\|\int_0^t\rme^{(t-\tau)\Delta}(\partial_1-\partial_2)(\rme^{\tau\Delta}u_1^0\rme^{\tau\Delta}u_1^0)\rmd\tau
  \Big\|_{B^{-1,\sigma}_{\infty\,q}(\mathscr{B}_m)}
\nonumber\\
  &\lesssim\; &m^{-2\sigma-\frac{2}{q}}\Big[\sum_{j\in\mathscr{B}_m}\Big(2^{-j}j^{\sigma}\Big\|\int_0^t\rme^{(t-\tau)\Delta}
  \partial_2V_{1j}\rmd\tau\Big\|_{\infty}\Big)^q\;\Big]^{\frac{1}{q}}
\nonumber\\
  &&\quad +m^{-2\sigma-\frac{2}{q}}\Big[\sum_{j\in\mathscr{B}_m}\Big(2^{-j}j^{\sigma}\Big\|\int_0^t\rme^{(t-\tau)\Delta}
  \partial_2V_{2j}\rmd\tau\Big\|_{\infty}\Big)^q\;\Big]^{\frac{1}{q}}
\nonumber\\ [0.2cm]
  &:=\; &I\!\!I\!\!I+I\!V.
\end{eqnarray}
  As for $\displaystyle\int_0^t\rme^{(t-\tau)\Delta}(\partial_1-\partial_2)U_{1j}\rmd\tau$ we have
\begin{eqnarray*}
  \int_0^t\rme^{(t-\tau)\Delta}\partial_2V_{1j}\rmd\tau&\;=\;&
  \bfi\mathscr{F}^{-1}\Big[\sum_{k\in\mathscr{A}_m}2^{2k+1}\int_0^t\rme^{-(t-\tau)|\xi|^2}\xi_2\Big(\rme^{-\tau|\xi|^2}
  \frac{\xi_2-\xi_1}{\xi_2}\Phi_{kj}^{++}\ast\rme^{-\tau|\xi|^2}\Phi_{kj}^{-+}
\\
  &&\quad +\rme^{-\tau|\xi|^2}\frac{\xi_2-\xi_1}{\xi_2}\Phi_{kj}^{+-}\ast
  \rme^{-\tau|\xi|^2}\Phi_{kj}^{--}\Big)\rmd\tau\Big]
\\
  &\;=\;&\bfi\mathscr{F}^{-1}\Big[\rme^{\bfi c_j\xi}\rme^{-t|\xi|^2}\xi_2\Big(\int_{\mathbb{R}^n}H^+_{j}(t,\xi,\eta)
  \rho(\xi-\eta-2b_j)\rho(\eta)\rmd\eta
\\
  &&\quad +\int_{\mathbb{R}^n}H^-_{j}(t,\xi,\eta)\rho(\xi-\eta+2b_j)\rho(\eta)\rmd\eta\Big)\Big],
\end{eqnarray*}
  where
\begin{eqnarray*}
  H^+_{j}(t,\xi,\eta)&\;:=\;&\sum_{k\in\mathscr{A}_m}2^{2k+1}
  \frac{\xi_2-\eta_2-\xi_1+\eta_1-2^{j-1}\varepsilon}{\xi_2-\eta_2+2^kn^{-\frac{1}{2}}-2^j\varepsilon}
  \int_0^t\rme^{-\tau(|\xi-\eta+a_k-b_j|^2+|\eta-a_k+b_j|^2-|\xi|^2)}\rmd\tau,
\\
  H^-_{j}(t,\xi,\eta)&\;:=\;&\sum_{k\in\mathscr{A}_m}2^{2k+1}
  \frac{\xi_2-\eta_2-\xi_1+\eta_1+2^{j-1}\varepsilon}{\xi_2-\eta_2+2^kn^{-\frac{1}{2}}+2^j\varepsilon}
  \int_0^t\rme^{-\tau(|\xi-\eta+a_k+b_j|^2+|\eta-a_k-b_j|^2-|\xi|^2)}\rmd\tau.
\end{eqnarray*}
  It is easy to see that on the supports of $\rho(\xi-\eta\mp2b_j)\rho(\eta)$ there respectively hold
$$
  \Big|\frac{\xi_2-\eta_2-\xi_1+\eta_1\mp2^{j-1}\varepsilon}{\xi_2-\eta_2+2^kn^{-\frac{1}{2}}\mp2^j\varepsilon}\Big|
  \lesssim 2^{j-k}\varepsilon,
$$
  so that similarly as before we have
$$
  |H^{\pm}_{j}(t,\xi,\eta)|\lesssim \sum_{k\in\mathscr{A}_m}2^{2k+1}\cdot 2^{j-k}\varepsilon\cdot
  \frac{1-\rme^{-t2^{2k+1}}}{2^{2k+1}}\lesssim 2^{j-16m}\varepsilon.
$$
  Hence
\begin{eqnarray*}
  \Big\|\int_0^t\rme^{(t-\tau)\Delta}\partial_2V_{1j}\rmd\tau\Big\|_{\infty}
  &\;\lesssim\;&\int_{\mathbb{R}^n}\rme^{-t|\xi|^2}\xi_2\Big(\int_{\mathbb{R}^n}|H^+_{j}(t,\xi,\eta)|
  \rho(\xi\!-\!\eta\!-\!2b_j)\rho(\eta)\rmd\eta
\\
  &&\quad +\int_{\mathbb{R}^n}|H^-_{j}(t,\xi,\eta)|\rho(\xi\!-\!\eta\!+\!2b_j)\rho(\eta)\rmd\eta\Big)\rmd\xi
\\
  &\;\lesssim\;& 2^{2j-16m}\varepsilon^2, \quad \forall j\in\mathscr{B}_m.
\end{eqnarray*}
  It follows that
\begin{equation}
  I\!\!I\!\!I\lesssim m^{-2\sigma-\frac{2}{q}}\cdot m^{\sigma}2^{-16m}\varepsilon^2\cdot\Big(\sum_{j\in\mathscr{B}_m}2^{jq}\Big)^{\frac{1}{q}}
  \lesssim \varepsilon^22^{-8m}m^{-\sigma-\frac{2}{q}}.
\end{equation}
  For $\displaystyle\int_0^t\rme^{(t-\tau)\Delta}\partial_2V_{2j}\rmd\tau$ we have
\begin{eqnarray*}
  \int_0^t\rme^{(t-\tau)\Delta}\partial_2V_{2j}\rmd\tau&\;=&\sum_{{j'\in\mathscr{B}_m\atop j'<j}}\!\theta_{kjj'}
  \bfi\mathscr{F}^{-1}\Big[\rme^{\bfi c_j\xi}\rme^{-t|\xi|^2}\xi_2\Big(\int_{\mathbb{R}^n}\!\rme^{\bfi(c_j-c_{j'})\eta}H^+_{j'}(t,\xi,\eta)
  \rho(\xi\!-\!\eta\!-\!b_j\!-\!b_{j'})\rho(\eta)\rmd\eta
\\
  &&\quad +\int_{\mathbb{R}^n}\!\rme^{\bfi(c_j-c_{j'})\eta}H^-_{j'}(t,\xi,\eta)\rho(\xi\!-\!\eta\!+\!b_j\!+\!b_{j'})\rho(\eta)\rmd\eta\Big)\Big].
\end{eqnarray*}
  It follows that
\begin{eqnarray*}
  \Big\|\int_0^t\rme^{(t-\tau)\Delta}\partial_2V_{2j}\rmd\tau\Big\|_{\infty}
  &\;\lesssim\;&\sum_{{j'\in\mathscr{B}_m\atop j'<j}}\!2^j\varepsilon\cdot 2^{j'-16m}\varepsilon
  \lesssim 2^{2j-16m}\varepsilon^2, \quad \forall j\in\mathscr{B}_m.
\end{eqnarray*}
  Consequently,
\begin{equation}
  I\!V\lesssim m^{-2\sigma-\frac{2}{q}}\cdot m^{\sigma}2^{-16m}\varepsilon^2\cdot\Big(\sum_{j\in\mathscr{B}_m}2^{jq}\Big)^{\frac{1}{q}}
  \lesssim \varepsilon^22^{-8m}m^{-\sigma-\frac{2}{q}}.
\end{equation}
  Substituting (4.13) and (4.14) into (4.12), we get (4.11). $\quad\Box$
\medskip

  {\bf Lemma 4.4} \ \ {\em Let $m\gg-\ln\varepsilon$ and $t=\varepsilon 2^{-32m}$. Then for any $1\leqslant q\leqslant\infty$ and $t>0$ we have}
\begin{equation}
  \Big\|\int_0^t\rme^{(t-\tau)\Delta}\partial_1\!\sum_{\alpha,\beta=1}^2\!\frac{\partial_{\alpha}\partial_{\beta}}{\Delta}
  (\rme^{\tau\Delta}u_{\alpha}^0\rme^{\tau\Delta}u_{\beta}^0)
  \rmd\tau\Big\|_{B^{-1,\sigma}_{\infty\,q}(\mathscr{B}_m)}\lesssim\varepsilon^3 m^{1-\sigma-\frac{1}{q}}.
\end{equation}

  {\em Proof}:\ \ A basic observation in getting (4.15) is that on the supports of $\rho(\xi-\eta\mp 2b_j)\rho(\eta)$ and
  $\rho(\xi-\eta\mp b_j\mp b_{j'})\rho(\eta)$, the symbol of the pseudo-differential operator $\displaystyle
  \frac{\partial_{\alpha}\partial_{\beta}}{\Delta}$ is bounded by $C\varepsilon^2$, where $C$ represents positive constant. Indeed, on these
  supports we have
$$
  \xi_\alpha\sim 2^j\varepsilon \; (\alpha=1,2)\;\; \mbox{and} \;\; |\xi|\sim 2^j\;\; (\mbox{assuming $j'<j$}), \quad \mbox{so that} \;\;
  \Big|\frac{\xi_{\alpha}\xi_{\beta}}{|\xi|^2}\Big|\lesssim\varepsilon^2 \; (\alpha,\beta=1,2).
$$
  There are four terms in the sum $\displaystyle\sum_{\alpha,\beta=1}^2$, and consequently the left-hand side of (4.15) can be bounded by a sum of
  four terms. We estimate each term separately.

  (1)\ \ Estimate of the term with $\alpha=\beta=1$. This term can be estimated as in the proof of (4.5). Indeed, similarly as in (4.6) we have
$$
  \Big\|\int_0^t\rme^{(t-\tau)\Delta}\partial_1\frac{\partial_1\partial_1}{\Delta}
  (\rme^{\tau\Delta}u_1^0\rme^{\tau\Delta}u_1^0)\rmd\tau\Big\|_{B^{-1,\sigma}_{\infty\,q}(\mathscr{B}_m)}
  \lesssim V+V\!I,
$$
  where $V$ and $V\!I$ are expressions obtained from modifying $I$ and $I\!\!I$, respectively, by replacing $\partial_1-\partial_2$ with
  $\displaystyle\partial_1\frac{\partial_{\alpha}\partial_{\beta}}{\Delta}$. By using some similar argument as in the proof of (4.9) (and using
  the inequalities $1-\rme^{-t2^{32m}}\leqslant 1$ and $\rme^{-t2^{2j}}\leqslant 1$) we have
$$
  V\lesssim m^{-2\sigma-\frac{2}{q}}\cdot\varepsilon^3 m^{1+\sigma+\frac{1}{q}}\lesssim\varepsilon^3m^{1-\sigma-\frac{1}{q}},
$$
  and using some similar argument as in the proof of (4.10) we have
$$
  V\!I\lesssim m^{-2\sigma-\frac{2}{q}}\cdot 2^{-4m}m^{\sigma}\cdot\varepsilon^3 m
  \lesssim\varepsilon^3 m^{1-\sigma-\frac{2}{q}}2^{-4m}.
$$
  Hence
\begin{equation}
  \Big\|\int_0^t\rme^{(t-\tau)\Delta}\partial_1\frac{\partial_1\partial_1}{\Delta}
  (\rme^{\tau\Delta}u_1^0\rme^{\tau\Delta}u_1^0)\rmd\tau\Big\|_{B^{-1,\sigma}_{\infty\,q}(\mathscr{B}_m)}
  \lesssim\varepsilon^3 m^{1-\sigma-\frac{1}{q}}.
\end{equation}

  (2)\ \ Estimate of the term with $\alpha=1$, $\beta=2$ or $\alpha=2$, $\beta=1$. We have
$$
  \rme^{\tau\Delta}u_1^0\rme^{\tau\Delta}u_2^0=\rme^{\tau\Delta}u_1^0\rme^{\tau\Delta}(u_1^0+u_2^0)-\rme^{\tau\Delta}u_1^0\rme^{\tau\Delta}u_1^0.
$$
  By using some similar argument as in the proof of (4.11) we have
$$
  \Big\|\int_0^t\rme^{(t-\tau)\Delta}\partial_1\frac{\partial_1\partial_2}{\Delta}
  [\rme^{\tau\Delta}u_1^0\rme^{\tau\Delta}(u_1^0+u_2^0)]\rmd\tau\Big\|_{B^{-1,\sigma}_{\infty\,q}(\mathscr{B}_m)}
  \lesssim\varepsilon^4 2^{-8m}m^{-\sigma-\frac{2}{q}},
$$
  and similar to (4.16) we have
$$
  \Big\|\int_0^t\rme^{(t-\tau)\Delta}\partial_1\frac{\partial_1\partial_2}{\Delta}
  (\rme^{\tau\Delta}u_1^0\rme^{\tau\Delta}u_1^0)\rmd\tau\Big\|_{B^{-1,\sigma}_{\infty\,q}(\mathscr{B}_m)}
  \lesssim\varepsilon^3 m^{1-\sigma-\frac{1}{q}}.
$$
  Hence
\begin{equation}
  \Big\|\int_0^t\rme^{(t-\tau)\Delta}\partial_1\frac{\partial_1\partial_2}{\Delta}
  (\rme^{\tau\Delta}u_1^0\rme^{\tau\Delta}u_2^0)\rmd\tau\Big\|_{B^{-1,\sigma}_{\infty\,q}(\mathscr{B}_m)}
  \lesssim\varepsilon^3 m^{1-\sigma-\frac{1}{q}}.
\end{equation}

  (3)\ \ Estimate of the term with $\alpha=\beta=2$. We have
$$
  \rme^{\tau\Delta}u_2^0\rme^{\tau\Delta}u_2^0=\rme^{\tau\Delta}(u_1^0+u_2^0)\rme^{\tau\Delta}(u_1^0+u_2^0)
  -2\rme^{\tau\Delta}u_1^0\rme^{\tau\Delta}(u_1^0+u_2^0)+\rme^{\tau\Delta}u_1^0\rme^{\tau\Delta}u_1^0.
$$
  By using some similar argument as in the proof of (4.11) we have
$$
  \Big\|\int_0^t\rme^{(t-\tau)\Delta}\partial_1\frac{\partial_1\partial_2}{\Delta}
  [\rme^{\tau\Delta}(u_1^0+u_2^0)\rme^{\tau\Delta}(u_1^0+u_2^0)]\rmd\tau\Big\|_{B^{-1,\sigma}_{\infty\,q}(\mathscr{B}_m)}
  \lesssim\varepsilon^5 2^{-16m}m^{-\sigma-\frac{2}{q}}.
$$
  Hence
\begin{eqnarray}
  &&\Big\|\int_0^t\rme^{(t-\tau)\Delta}\partial_1\frac{\partial_2\partial_2}{\Delta}
  (\rme^{\tau\Delta}u_2^0\rme^{\tau\Delta}u_2^0)\rmd\tau\Big\|_{B^{-1,\sigma}_{\infty\,q}(\mathscr{B}_m)}
\nonumber\\
  &\;\lesssim\;&\varepsilon^5 2^{-16m}m^{-\sigma-\frac{2}{q}}+\varepsilon^4 2^{-8m}m^{-\sigma-\frac{2}{q}}
  +\varepsilon^3 m^{1-\sigma-\frac{1}{q}}
  \lesssim\varepsilon^3 m^{1-\sigma-\frac{1}{q}}.
\end{eqnarray}

  Now, summing up (4.16), (4.17) and (4.18) (with (4.17) twice), we obtain (4.15). $\quad\Box$
\medskip

  We are now ready to give the proof of Theorem 1.3.
\medskip

  {\em Proof of Theorem 1.3}:\ \ As we mentioned in the beginning of this section, here we only consider the case $2<q\leqslant\infty$,
  $1-2/q\leqslant\sigma<1-1/q$. Hence, in what follows we assume that these conditions are satisfied.

  Given $T>0$, let $X$ be the following function space on $\mathbb{R}^n\times [0,T]$:
$$
  X=\{u\in L^{\infty}_{\rm loc}((0,T],L^{\infty}(\mathbb{R}^n)):\; \|u\|_X<\infty\},
$$
  where
$$
  \|u\|_X:=\sup_{0<t<T}\sqrt{t}\|u(\cdot,t)\|_{L^{\infty}(\mathbb{R}^n)}+\sup_{{x\in\mathbb{R}^n\atop 0<R<\sqrt{T}}}
  \Big(\frac{1}{|B(x,R)|}\int_0^{R^2}\!\!\!\int_{B(x,R)}|u(y,t)|^2\rmd y\rmd t\Big)^{\frac{1}{2}},
$$
  and let $Y=X\cap L^{\infty}((0,T),bmo^{-1})$ with standard norm for joint space. It is well-known that there exists $\epsilon_T>0$ such that for
  $\bfu_0\in bmo^{-1}$ with $\|\delta\bfu_0\|_{bmo^{-1}}\leqslant \epsilon_T$, the problem (4.1) has a unique solution $\bfu=\bfu(\delta,t)\in Y$, and
\begin{equation}
  \|\bfu(\delta,\cdot)\|_Y\lesssim\delta\|\bfu_0\|_{bmo^{-1}}, \quad \forall\delta\in (0,\delta_0).
\end{equation}
  Moreover, there holds the following estimate:
\begin{equation}
  \|B(\bfu,\bfv)\|_Y\lesssim\|\bfu\|_X\|\bfv\|_X, \quad \forall\bfu,\bfv\in X.
\end{equation}
  We refer the reader to see Theorem 16.1, Lemma 16.3 and the proof of the corollary following Theorem 16.2 of \cite{LEM02} for proofs of these
  assertions. Since $\bfu=\delta\rme^{t\Delta}\bfu_0+B(\bfu,\bfu)$, from (4.19) and (4.20) it follows that
\begin{equation}
  \|\bfu(\delta,\cdot)-\delta\rme^{t\Delta}\bfu_0\|_Y\lesssim\|\bfu(\delta,\cdot)\|_X^2\lesssim\delta^2\|\bfu_0\|_{bmo^{-1}}^2,
  \quad \forall\delta\in (0,\delta_0).
\end{equation}
  Let $\bfv(\delta,t)=B(\rme^{t\Delta}\bfu_0,\rme^{t\Delta}\bfu_0)$ and set
\begin{eqnarray*}
  \bfw(\delta,t)&\;=\;&\bfu(\delta,t)-\delta\rme^{t\Delta}\bfu_0-\delta^2\bfv(\delta,t)
\\
  &\;=\;&\int_0^t\rme^{(t-\tau)\Delta}\mathbb{P}\nabla\cdot[\bfu(\delta,\tau)\otimes\bfu(\delta,\tau)
  -\delta^2\rme^{\tau\Delta}\bfu_0\otimes\rme^{\tau\Delta}\bfu_0]\rmd\tau.
\end{eqnarray*}
  Since
$$
  \bfw(\delta,t)=\int_0^t\rme^{(t-\tau)\Delta}\mathbb{P}\nabla\cdot\{[\bfu(\delta,\tau)-\delta\rme^{t\Delta}\bfu_0]\otimes\bfu(\delta,\tau)
  +\delta\rme^{t\Delta}\bfu_0\otimes[\bfu(\delta,\tau)-\delta\rme^{t\Delta}\bfu_0]\}\rmd\tau,
$$
  by (4.19) $\sim$ (4.21) it follows that
\begin{equation}
  \|\bfw(\delta,\cdot)\|_Y\lesssim(\|\bfu(\delta,\cdot)\|_X+\delta\|\rme^{t\Delta}\bfu_0\|_X)
  \|\bfu(\delta,\cdot)-\delta\rme^{t\Delta}\bfu_0\|_X\lesssim\delta^3\|\bfu_0\|_{bmo^{-1}}^3.
\end{equation}
  Now let $\bfu_0$ be the vector function given by (4.2). Since $\bfu_0\in S(\mathbb{R}^n))$, the above results apply to it. By applying Lemma 4.1 we
  have
\begin{equation}
  \|\bfu_0\|_{bmo^{-1}}\lesssim\|\bfu_0\|_{B^{-1}_{\infty\,2}}\lesssim m^{\frac{1}{2}-\frac{1}{q}-\sigma}\|\bfu_0\|_{B^{-1,\sigma}_{\infty\,q}}
  \lesssim m^{\frac{1}{2}-\frac{1}{q}-\sigma}.
\end{equation}
  Thus by (4.22) we have
\begin{equation}
  \|\bfw(\delta,\cdot)\|_Y\lesssim\delta^3 m^{\frac{3}{2}-\frac{3}{q}-3\sigma},
\end{equation}
  which implies that
\begin{equation}
  \|\bfw(\delta,t)\|_{B^{-1,\sigma}_{\infty\,q}(\mathscr{B}_m)}\lesssim m^{\sigma+\frac{1}{q}}
  \|\bfw(\delta,t)\|_{B^{-1}_{\infty\infty}(\mathscr{B}_m)}\lesssim m^{\sigma+\frac{1}{q}}\|\bfw(\delta,\cdot)\|_Y
  \lesssim\delta^3 m^{\frac{3}{2}-\frac{2}{q}-2\sigma}, \quad \forall t\in (0,T).
\end{equation}
  Let $v_1(\delta,t)$ be the first component of $\bfv(\delta,t)$. From Lemmas 4.2 $\sim$ 4.4 we see that with $t=\varepsilon 2^{-32m}$ and
  $m\gg-\ln\varepsilon$,
$$
  \|\bfv(\delta,t)\|_{B^{-1,\sigma}_{\infty\,q}(\mathscr{B}_m)}\gtrsim\|v_1(\delta,t)\|_{B^{-1,\sigma}_{\infty\,q}(\mathscr{B}_m)}
  \gtrsim\varepsilon^2 m^{1-\sigma-\frac{1}{q}}.
$$
  Hence
\begin{eqnarray}
  \|\bfu(\delta,t)\|_{B^{-1,\sigma}_{\infty\,q}}&\;\gtrsim\;&\|\bfw(\delta,t)-\delta^2\bfv(\delta,t)\|_{B^{-1,\sigma}_{\infty\,q}(\mathscr{B}_m)}
  -\delta\|\rme^{t\Delta}\bfu_0\|_{B^{-1,\sigma}_{\infty\,q}}
\nonumber\\
  &\;\gtrsim\;&\delta^2\|\bfv(\delta,t)\|_{B^{-1,\sigma}_{\infty\,q}(\mathscr{B}_m)}-\|\bfw(\delta,t)\|_{B^{-1,\sigma}_{\infty\,q}(\mathscr{B}_m)}
  -\delta\|\bfu_0\|_{B^{-1,\sigma}_{\infty\,q}}
\nonumber\\
  &\;\geqslant\;& C_0\delta^2\varepsilon^2 m^{1-\sigma-\frac{1}{q}}-C_1\delta^3 m^{\frac{3}{2}-\frac{2}{q}-2\sigma}-C_2\delta
  \quad \mbox{for}\;\; t=\varepsilon 2^{-32m}.
\end{eqnarray}
  Since the conditions $1-2/q\leqslant\sigma<1-1/q$ and $q>2$ imply that $1-\sigma-1/q>0$ and $1-\sigma-1/q>3/2-2/q-2\sigma$, it follows that for any
  given $0<\delta<\delta_0$ and $0<\varepsilon\ll 1$, by choosing $m$ sufficiently large (depending on $\delta$ and $\varepsilon$), we have
$$
  \|\bfu(\delta,t)\|_{B^{-1,\sigma}_{\infty\,q}}\geqslant\frac{1}{2}C_0\delta^2\varepsilon^2 m^{1-\sigma-\frac{1}{q}}
  \quad \mbox{for}\;\; t=\varepsilon 2^{-32m}.
$$
  This completes the proof. $\quad\Box$
\medskip

  {\bf Remark 4.5} \ \ Note that the above proof actually works under the weaker conditions $2<q\leqslant\infty$ and $1/2-1/q\leqslant\sigma<1-1/q$.
\medskip

  {\bf Remark 4.6} \ \ To treat the case $1\leqslant q\leqslant 2$ and $0\leqslant\sigma<1/q$, we need to use a different class of initial values
  $\bfu_0$ which are obtained by modifying the definition of (4.2) as follows: In the first two lines of (4.2), remove the first sum $\displaystyle
  \sum_{k\in\mathscr{A}_m}$, replace $m^{-\sigma-\frac{1}{q}}$ with $m^{-\sigma}$, and put $k=16m$. The arguments in the statements and proofs of
  Lemmas 4.1 $\sim$ 4.4 must be correspondingly modified: The right-hand sides of (4.5), (4.11) and (4.15) need be replaced with $\varepsilon^2
  m^{\frac{1}{q}-\sigma}$, $\varepsilon^2 2^{-8m}m^{-\sigma}$ and $\varepsilon^4 m^{\frac{1}{q}-\sigma}$, respectively. We omit the details here.
  Note that this can also be regarded as a modification to the argument of \cite{Wang}.
\medskip

  {\bf Remark 4.7} \ \ To treat the case $2<q\leqslant\infty$ and $0\leqslant\sigma<1-2/q$, another different class of initial values $\bfu_0$ have
  to be employed, which are obtained by modifying the definition of (4.2) in another way as follows: In the first two lines of (4.2), remove the
  second sum $\displaystyle\sum_{l\in\mathscr{B}_m}$, and put $l=4m$. The corresponding modifications in the statements of Lemmas 4.2
  $\sim$ 4.4 are as follows: All $\mathscr{B}_m$ in these lemmas need be replaced with the single-point set $\{4m\}$, and the right-hand sides of
  (4.5), (4.11) and (4.15) need be replaced with $\varepsilon^2m^{1-\sigma-\frac{2}{q}}$, $\varepsilon^2 2^{-12m}m^{-\sigma-\frac{2}{q}}$ and
  $\varepsilon^3 m^{1-\sigma-\frac{2}{q}}$, respectively. In this case, the estimates in (4.23) $\sim$ (4.26) need also be correspondingly modified:
  Not only the set $\mathscr{B}_m$ need be replaced with the single-point set $\{4m\}$, but also the final bounds in these estimates need be replaced
  with $m^{-\sigma}$, $\delta^3m^{-3\sigma}$, $\delta^3m^{-2\sigma}$ and
$$
  C_0\delta^2\varepsilon^2 m^{1-\sigma-\frac{2}{q}}-C_1\delta^3 m^{-2\sigma}-C_2\delta,
$$
  respectively. We omit the details here.
\medskip

  {\bf Remark 4.8} \ \ For the case $2<q\leqslant\infty$ and $0\leqslant\sigma<1-2/q$, an alternative proof is to modify the argument of
  \cite{Yon10}. The modification is as follows: Let
\begin{eqnarray*}
  & v^0=(0,0,1,0,\cdots,0), \quad  v^1=(0,1,0,0,\cdots,0), \quad w^0=(1,0,0,0,\cdots,0), \quad  w^1=(0,0,1,0,\cdots,0), &
\\
  & a_k^0=2^{4(m+k)}w^0, \quad  a_k^1=2^{4(m+k)}w^0+2^{m}w^1,  \quad k=1,2,\cdots,m.
\end{eqnarray*}
  Then set
$$
  \bfu_0(x)=m^{-\sigma-\frac{1}{q}}\sum_{k=1}^m|a_k^0|[v^0\cos(a_k^0x)+v^1\cos(a_k^1x)].
$$
  Using this function as the initial value and correspondingly making necessary modifications to the argument of \cite{Yon10}, we obtain a different
  proof to ill-posedness of the problem (1.2) in $B^{-1,\sigma}_{\infty\, q}(\mathbb{R}^n)$ for the case $2<q\leqslant\infty$ and
  $0\leqslant\sigma<1-2/q$.
\medskip

  {\bf Remark 4.9} \ \ Here we only treated the case $n\geqslant 3$. After making some modifications to the argument given above as in Section 3 of
  \cite{Wang}, we see that Theorem 1.3 also holds for the case $n=2$.


{\small
\begin{thebibliography}{}
\bibitem{BAR96} O. Barraza, Self-similar solutions in weak $L^p$-spaces of the Navier-Stokes equations, \textit{Revista Mat Iberoamer.},
 \textbf{12}(1996), 411--439.
\bibitem{BerLof} J. Bergh and J. L\"{o}fstr\"{o}m, \textit{Interpolation Spaces: An Introduction}, Springer--Verlag, New York, 1976.
\bibitem{BP08} J. Bourgain and J. N. Pavlovi\'{c}, Ill-posedness of the Navier-Stokes equations in a critical space in 3D, \textit{J. Funct. Anal.},
 \textbf{255}(2008), 2233--2247.
\bibitem{CAN97} M. Cannone, A generalization of a theorem by Kato on Navier-Stokes equations, \textit{Rev. Mat. Iberoam}, \textbf{13}(1997), 515--541.
\bibitem{CAN04} M. Cannone, Harmonic analysis tools for solving the incompressible Navier-Stokes equations, in \textit{Handbook of Mathematical Fluid
  Dynamics vol. III} (S. Friedlander and D. Serre edit.), pp. 161--244, Elsevier, north Holland:  2004.
\bibitem{Cui} S. Cui, Weak solutions for the Navier-Stokes equations with $B^{-1,\sigma}_{\infty\, q}+B_{X_r}^{-1+r,\frac{2}{1-r}}+L^2$
  initial data, \textit{J. Math. Phys.}, \textbf{54}(2013), 051503-1--18.
\bibitem{FABJR72} E, Fabes, B. Johns and n. Riviere, The initial value problem for the Navier-Stokes equations with data in $L^p$, \textit{Arch. Rat.
  Mech. Anal.}, \textbf{45}(1972), 222--240.
\bibitem{FUJK64} H. Fujita and T. Kato, On the Navier-Stokes initial value problem I, \textit{Arch. Ration. Mech. Anal.}, \textbf{16}(1964), 269--315.
\bibitem{GIG86} Y. Giga, Solutions of semilinear parabolic equations in $L^p$ and regularity of weak solutions of the Navier-Stokes system, \textit{J.
  Diff. Equa.}, \textbf{62}(1986), 186--212.
\bibitem{GM89} Y. Giga and T. Miyakawa, Navier-Stokes flows in $\mathbb{R}^3$ with measurea as initial vorticity and the Morrey spaces, \textit{Comm.
  P. D. E.}, \textbf{14}(1989), 577--618.
\bibitem{KAT84} T. Kato, Strong $L^p$ solutions of the Navier-Stokes equations in $\mathbb{R}^m$ with applications to weak solutions, \textit{Math. Z.},
  \textbf{187}(1984), 471--480.
\bibitem{KOCT01}  H. Koch and D. Tataru, Well-posedness for the Navier-Stokes equations, \textit{Adv. in Math.}, \textbf{157}(2001), 22--35.
\bibitem{LEM02} P. G. Lemari\'{e}-Rieusset, \textit{Recent developments in the Navier-Stokes problems}, Research notes in Mathematics, Chapman \&
  Hall/CRC, 2002.
\bibitem{LEM07} P. G. Lemari\'{e}-Rieusset, The Navier-Stokes equations in the critical Morrey-Campanato space, \textit{Revista Mat. Iberoamer},
  \textbf{23}(2007), 1--31.
\bibitem{LEM16} P. G. Lemari\'{e}-Rieusset, \textit{The Navier-Stokes Problem in the 21st Century}, Boca Raton: Taylor \& Francis Group, 2016.
\bibitem{PLA96} F. Planchon, Solutions globales et comportement asymptotique pour les $\acute{e}$quations de Navier-Stokes, Th$\grave{e}$se, Ecole
  Polytechnique, 1996.
\bibitem{TER99} E. Terraneo, Application de certains espaces de l'analyse harmonique r$\acute{e}$elle $\grave{a}$ l'$\acute{e}$tude des
  $\acute{e}$quations de Navier-Stokes et de la chaleur nonlin$\acute{e}$aires, Th$\grave{e}$se, Univ. Evry., 1999.
\bibitem{Wang} B. Wang, Ill-posedness for the Navier-Stokes equations in critical Besov spaces $B^{s}_{\infty,q}$, \textit{Advances in Math.},
  \textbf{268}(2015), 350--372.
\bibitem{WEI81} F. B. Weissler, The Navier-Stokes initial value problem in $L^p$, \textit{Arch. Rational Mech. Anal.}, \textbf{74}(1981), 219--230.
\bibitem{Yon10} T. Yoneda, Ill-posedness of the 3D-Navier-Stokes equations in a generalized Besov space near $BMO^{-1}$, \textit{J. Func. Anal.},
  \textbf{258}(2010), 3376--3387.
\end{thebibliography}
}
\end{document}